\documentclass[a4paper,12pt]{article}
\usepackage[utf8]{inputenc}
\usepackage[T1]{fontenc}
\usepackage{amsfonts,amsmath,amssymb,amsthm}
\theoremstyle{definition}

\theoremstyle{plain}

\newtheorem{lemma}{Lemma}
\newtheorem{proposition}{Proposition}
\newtheorem{theorem}{Theorem}
\usepackage{graphicx}
\usepackage{caption}
\usepackage{pgf,tikz,circuitikz,subfig}
\usetikzlibrary{arrows,shapes,positioning}
\usepackage{booktabs,multirow,tabularx,tablefootnote}
\usepackage{hyperref}
\hypersetup{colorlinks,%
	citecolor={red}, 
	urlcolor={blue},
	linkcolor={blue},
	breaklinks={true},
	pagebackref={true},
	hyperindex={true},
}
\usepackage{url}
\usepackage{cite}
\usepackage[left=2cm,right=2cm,top=2cm,bottom=4cm]{geometry}
\usepackage[ruled]{algorithm}
\usepackage{algorithmic}
\usepackage{newtxtext}
\usepackage{courier}
\usepackage{enumitem}
\newlist{abbrv}{itemize}{1}
\setlist[abbrv,1]{label=,labelwidth=0.9in,align=parleft,noitemsep,leftmargin=!}

\newcommand{\rv}[1]{\boldsymbol{#1}}

\newcommand{\ub}[1]{\overline{#1}}

\newcommand{\geo}[1]{\mathtt{#1}}
\DeclareMathOperator{\subj}{s.t.}

\title{Maximal perimeter and maximal width of a convex small polygon}
\author{Christian Bingane\thanks{Department of Mathematics and Industrial Engineering, Polytechnique Montreal, Montreal, Quebec, Canada, H3C 3A7. Email: \url{christian.bingane@polymtl.ca}}}
\begin{document}
\maketitle
\begin{abstract}
A small polygon is a polygon of unit diameter. The maximal perimeter and the maximal width of a convex small polygon with $n=2^s$ sides are unknown when $s \ge 4$. In this paper, we propose an approach to construct convex small $n$-gons of large perimeter and large width when $n=2^s$ with $s\ge 2$. Assuming the existence of an axis of symmetry, a convex small $n$-gon is described as a composition of $n/2$ and both its perimeter and its width are given as functions of a single variable. By selecting the composition that minimizes the violation of a cycle constraint by a particular solution, the $n$-gons constructed outperform the best $n$-gons found in the literature. For example, for $n=64$, the perimeter and the width obtained are within $10^{-22}$ and $10^{-12}$ of the maximal perimeter and the maximal width, respectively. From our results, it appears that Mossinghoff's conjecture on the diameter graph of a convex small $2^s$-gon with maximal perimeter is not true when $s \ge 4$.
\end{abstract}


\section{Introduction}
Let $\geo{P}$ be a convex polygon. The {\em diameter} of $\geo{P}$ is the maximum distance between pairs of its vertices. The polygon $\geo{P}$ is {\em small} if its diameter equals one. The diameter graph of a small polygon is defined as the graph with the vertices of the polygon, and an edge between two vertices exists only if the distance between these vertices equals one. Diameter graphs of some convex small polygons are represented in Figure~\ref{figure:4gon}, Figure~\ref{figure:6gon}, and Figure~\ref{figure:8gon}. The solid lines illustrate pairs of vertices which are unit distance apart. The height associated to a side of~$\geo{P}$ is defined as the maximum distance between a vertex of~$\geo{P}$ and the line containing the side. The minimum height for all sides is the {\em width} of the polygon~$\geo{P}$.

When $n=2^s$ with integer $s\ge 4$, both the maximal perimeter and the maximal width of a convex small $n$-gon are unknown. However, tight bounds may be obtained analytically. It is well known that, for an integer $n \ge 3$, the value $2n\sin \frac{\pi}{2n}$~\cite{reinhardt1922,vincze1950,datta1997} is an upper bound on the perimeter $L(\geo{P}_n)$ of a convex small $n$-gon $\geo{P}_n$ and the value $\cos \frac{\pi}{2n}$~\cite{gashkov1985,bezdek2000} an upper bound on its width $W(\geo{P}_n)$. Recently, the author~\cite{bingane2022e} constructed a family of convex small $n$-gons, for $n=2^s$ with $s\ge 3$, whose perimeters and widths differ from the upper bounds $2n\sin \frac{\pi}{2n}$ and $\cos \frac{\pi}{2n}$ by just $O(1/n^6)$ and $O(1/n^4)$, respectively. By constrast, both the perimeter and the width of a regular small $n$-gon differ by $O(1/n^2)$ when $n \ge 4$ is even. In the present paper, we further tighten lower bounds on the maximal perimeter and the maximal width. Following the strategy employed in~\cite{bingane2022e}, our first result is the following:

\begin{theorem}\label{thm:Cn}
	Suppose $n=2^s$ with integer $s\ge 4$. Let $\ub{L}_n := 2n \sin \frac{\pi}{2n}$ denote an upper bound on the perimeter $L(\geo{P}_n)$ of a convex small $n$-gon $\geo{P}_n$, and $\ub{W}_n := \cos \frac{\pi}{2n}$ denote an upper bound on its width $W(\geo{P}_n)$. Then there exists a convex small $n$-gon $\geo{C}_n$ such that
	\[
	\begin{aligned}
		\ub{L}_n - L(\geo{C}_n) &= \frac{\pi^9}{8n^8} + O\left(\frac{1}{n^{10}}\right),\\
		\ub{W}_n - W(\geo{C}_n) &= \frac{\pi^5}{4n^5} + O\left(\frac{1}{n^7}\right).
	\end{aligned}
	\]
\end{theorem}

For all $n = 2^s$ and $s \ge 4$, the diameter graph of the $n$-gon $\geo{C}_n$ has a cycle of length $3n/4-1$ plus $n/4+1$ pendant edges. In 2006, Mossinghoff~\cite{mossinghoff2006b} conjectured that, when $n = 2^s$ and $s \ge 3$, the diameter graph of a convex small $n$-gon of maximal perimeter has a cycle of length $n/2+1$ plus $n/2-1$ additional pendant edges, and that is verified for $s=3$. However, it appears that the conjecture is no longer true for $s\ge 4$ as the perimeter of $\geo{C}_n$ exceeds that of the optimal $n$-gon obtained by Mossinghoff.

One may ask if the diameter graph of a convex small $n$-gon of maximal perimeter or maximal width is the same as that of the $n$-gon $\geo{C}_n$ when $n=2^s$ and $s \ge 4$. It seems that it is not the case when $s\ge 5$. In 2021, Xue, Lian, Wang, and Zhang~\cite{xue2021} proposed a better family of convex of small $n$-gons, for $n=2^s \ge 4$, whose perimeters and widths are within $O(\pi^{2s}/n^{s+5})$ and $O(\pi^{s}/n^{\frac{s+7}{2}})$ of the maximal perimeter and the maximal width, respectively. The diameter graph of their $n$-gon has a cycle of length $\lfloor 2n/3\rceil$, which is a Jacobsthal number (sequence \href{https://oeis.org/A001045}{A001045} in the \href{https://oeis.org/}{OEIS}). However, we can still improve the lower bounds for $n\ge 32$.

Our second and main result is proposing an approach to construct a convex small $n$-gon of large perimeter and large width when $n=2^s$ with integer $s\ge 2$. This approach generalizes the strategy employed in~\cite{bingane2022e}. Assuming the existence of an axis of symmetry, we describe a convex small $n$-gon~$\geo{P}_n$ as a composition of $n/2$ and both its perimeter $L(\geo{P}_n)$ and its width $W(\geo{P}_n)$ are given as functions of a single variable. Then the convex small $n$-gon of large perimeter and large width, denoted $\geo{S}_n$, corresponds to the composition that minimizes the violation of a cycle constraint  by a particular solution among all $2^{n/2-1}-1$ compositions of $n/2$. For $n=64$, the perimeter and the width of $\geo{S}_{64}$ are within $10^{-22}$ and $10^{-12}$ of the maximal perimeter and the maximal width, respectively. In comparison, the perimeter and the width of the best prior $64$-gon, constructed in~\cite{xue2021}, are respectively within $10^{-14}$ and $10^{-8}$ of the maximal perimeter and the maximal width.

The remainder of this paper is organized as follows. Section~\ref{sec:ngon} recalls principal results on the maximal perimeter and the maximal width of convex small polygons. We develop our approach to construct convex small polygons of large perimeter and large width  in Section~\ref{sec:Cn} and the proof of Theorem~\ref{thm:Cn} is given. The convex small polygons of large perimeter and large width obtained are presented in Section~\ref{sec:Sn}. We maximize the perimeter and obtain polygons with even larger perimeters in Section~\ref{sec:maxperi}. We conclude the paper in Section~\ref{sec:conclusion}.

\begin{figure}[H]
	\centering
	\subfloat[$(\geo{R}_4,2.828427,0.707107)$]{
		\begin{tikzpicture}[scale=5]
			\draw[dashed] (0,0) -- (0.5000,0.5000) -- (0,1) -- (-0.5000,0.5000) -- cycle;
			\draw (0,0) -- (0,1);
			\draw (0.5000,0.5000) -- (-0.5000,0.5000);
		\end{tikzpicture}
	}
	\subfloat[$(\geo{T}_4,3.035276,0.866025)$]{
		\begin{tikzpicture}[scale=5]
			\draw[dashed] (0.5000,0.8660) -- (0,1) -- (-0.5000,0.8660);
			\draw (0,1) -- (0,0) -- (0.5000,0.8660) -- (-0.5000,0.8660) -- (0,0);
		\end{tikzpicture}
		\label{figure:4gon:T4}
	}
	\caption{Two convex small $4$-gons $(\geo{P}_4,L(\geo{P}_4),W(\geo{P}_4))$: (a) Regular $4$-gon; (b) Tamvakis $4$-gon~\cite{taylor1953, tamvakis1987}}
	\label{figure:4gon}
\end{figure}

\begin{figure}[H]
	\centering
	\subfloat[$(\geo{R}_6,3,0.866025)$]{
		\begin{tikzpicture}[scale=5]
			\draw[dashed] (0,0) -- (0.4330,0.2500) -- (0.4330,0.7500) -- (0,1) -- (-0.4330,0.7500) -- (-0.4330,0.2500) -- cycle;
			\draw (0,0) -- (0,1);
			\draw (0.4330,0.2500) -- (-0.4330,0.7500);
			\draw (0.4330,0.7500) -- (-0.4330,0.2500);
		\end{tikzpicture}
	}
	\subfloat[$(\geo{R}_{3,6},3.105829,0.965926)$]{
		\begin{tikzpicture}[scale=5]
			\draw[dashed] (0,0) -- (0.3660,0.3660) -- (0.5000,0.8660) -- (0,1) -- (-0.5000,0.8660) -- (-0.3660,0.3660) -- cycle;
			\draw (0,0) -- (0.5000,0.8660) -- (-0.5000,0.8660) -- cycle;
			\draw (0,0) -- (0,1);
			\draw (0.3660,0.3660) -- (-0.5000,0.8660);
			\draw (0.5000,0.8660) -- (-0.3660,0.3660);
		\end{tikzpicture}
		\label{figure:6gon:R36}
	}
	\caption{Two convex small $6$-gons $(\geo{P}_6,L(\geo{P}_6),W(\geo{P}_6))$: (a) Regular $6$-gon; (b) Reinhardt $6$-gon~\cite{reinhardt1922}}
	\label{figure:6gon}
\end{figure}

\begin{figure}[H]
	\centering
	\subfloat[$(\geo{R}_8,3.061467,0.923880)$]{
		\begin{tikzpicture}[scale=5]
			\draw[dashed] (0,0) -- (0.3536,0.1464) -- (0.5000,0.5000) -- (0.3536,0.8536) -- (0,1) -- (-0.3536,0.8536) -- (-0.5000,0.5000) -- (-0.3536,0.1464) -- cycle;
			\draw (0,0) -- (0,1);
			\draw (0.3536,0.1464) -- (-0.3536,0.8536);
			\draw (0.5000,0.5000) -- (-0.5000,0.5000);
			\draw (0.3536,0.8536) -- (-0.3536,0.1464);
		\end{tikzpicture}
	}
	\subfloat[$(\geo{B}_8,3.121062,0.977609)$]{
		\begin{tikzpicture}[scale=5]
			\draw[dashed] (0,0) -- (0.2957,0.2043) -- (0.5000,0.5000) -- (0.4114,0.9114) -- (0,1) -- (-0.4114,0.9114) -- (-0.5000,0.5000) -- (-0.2957,0.2043) -- cycle;
			\draw (0,0) -- (0.4114,0.9114) -- (-0.5000,0.5000) -- (0.5000,0.5000) -- (-0.4114,0.9114) -- cycle;
			\draw (0,0) -- (0,1);
			\draw (0.4114,0.9114) -- (-0.2957,0.2043);\draw (-0.4114,0.9114) -- (0.2957,0.2043);
		\end{tikzpicture}
		\label{figure:8gon:B8}
	}
	\subfloat[$(\geo{B}_8^*,3.121147,0.976410)$]{
		\begin{tikzpicture}[scale=5]
			\draw[dashed] (0,0) -- (0.2983,0.2128) -- (0.5000,0.5188) -- (0.4217,0.9067) -- (0,1) -- (-0.4217,0.9067) -- (-0.5000,0.5188) -- (-0.2983,0.2128) -- cycle;
			\draw (0,0) -- (0,1);
			\draw (0,0) -- (0.4217,0.9067) -- (-0.5000,0.5188) -- (0.5000,0.5188)-- (-0.4217,0.9067) -- cycle;
			\draw (0.4217,0.9067) -- (-0.2983,0.2128);\draw (-0.4217,0.9067) -- (0.2983,0.2128);
		\end{tikzpicture}
		\label{figure:8gon:B8*}
	}
	\caption{Three convex small $8$-gons $(\geo{P}_8,L(\geo{P}_8),W(\geo{P}_8))$: (a) Regular $8$-gon; (b) An $8$-gon of maximal width~\cite{audet2013}; (c) $8$-gon of maximal perimeter~\cite{griffiths1975, audet2007a}}
	\label{figure:8gon}
\end{figure}

\section{Perimeters and widths of convex small polygons}\label{sec:ngon}
Let $L(\geo{P})$ denote the perimeter of a polygon $\geo{P}$ and $W(\geo{P})$ its width. For a given integer $n\ge 3$, let $\geo{R}_n$ denote the regular small $n$-gon. We have
\[
L(\geo{R}_n) =
\begin{cases}
	2n\sin \frac{\pi}{2n} &\text{if $n$ is odd,}\\
	n\sin \frac{\pi}{n} &\text{if $n$ is even,}
\end{cases}
\]
and
\[
W(\geo{R}_n) =
\begin{cases}
	\cos \frac{\pi}{2n} &\text{if $n$ is odd,}\\
	\cos \frac{\pi}{n} &\text{if $n$ is even.}
\end{cases}
\]

When $n$ has an odd factor $m$, consider the family of convex equilateral small $n$-gons constructed as follows:
\begin{enumerate}
	\item Transform the regular small $m$-gon  $\geo{R}_m$ into a Reuleaux $m$-gon by replacing each edge by a circle's arc passing through its end vertices and centered at the opposite vertex;
	\item Add at regular intervals $n/m-1$ vertices within each arc;
	\item Take the convex hull of all vertices.
\end{enumerate}
These $n$-gons are denoted $\geo{R}_{m,n}$ and
\[
\begin{aligned}
	L(\geo{R}_{m,n}) &= 2n\sin \frac{\pi}{2n},\\
	W(\geo{R}_{m,n}) &= \cos \frac{\pi}{2n}.
\end{aligned}
\]
The $6$-gon $\geo{R}_{3,6}$ is illustrated in Figure~\ref{figure:6gon:R36}.

\begin{theorem}[Reinhardt~\cite{reinhardt1922}, Vincze~\cite{vincze1950}, Datta~\cite{datta1997}, Gashkov~\cite{gashkov1985}, Bezdek and Fodor \cite{bezdek2000}]\label{thm:perimeter}
	For all $n \ge 3$, let $L_n^*$ denote the maximal perimeter among all convex small $n$-gons, $W_n^*$ denote the maximal width, $\ub{L}_n := 2n \sin \frac{\pi}{2n}$ and $\ub{W}_n := \cos \frac{\pi}{2n}$.
	\begin{itemize}
		\item When $n$ has an odd factor $m$, $L_n^* = \ub{L}_n$ and $W_n^* = \ub{W}_n$ are simultaneously achieved by finitely many equilateral $n$-gons~\cite{gashkov2007,gashkov2013,mossinghoff2011,hare2013,hare2019}, including~$\geo{R}_{m,n}$. The optimal $n$-gon $\geo{R}_{m,n}$ is unique if $m$ is prime and $n/m \le 2$.
		\item When $n=4$, $L_4^* = 2+\sqrt{6}-\sqrt{2}$~ is only achieved by the non-equilateral $4$-gon $\geo{T}_4$~\cite{taylor1953,tamvakis1987}, represented in Figure~\ref{figure:4gon:T4}, and $W_4^* = \sqrt{3}/2$ is achieved by infinitely many $4$-gons, including~$\geo{T}_4$.
		\item When $n=8$, $L_8^* = 3.121147\ldots$ is only achieved by the non-equilateral $8$-gon $\geo{B}_8^*$~\cite{griffiths1975,audet2007a}, represented in Figure~\ref{figure:8gon:B8*}, and $W_8^* = 0.977608\ldots$ is achieved by infinitely many $8$-gons~\cite{audet2013}, including~$\geo{B}_8$ represented in Figure~\ref{figure:8gon:B8}.
		\item When $n=2^s$ with $s\ge 2$, $L(\geo{R}_n) < L_n^* < \ub{L}_n$ and $W(\geo{R}_n) < W_n^* < \ub{W}_n$.
	\end{itemize}
\end{theorem}

For $n=2^s$ with $s\ge 4$, tight lower bounds on the maximal perimeter and the maximal width may be obtained analytically. The author~\cite{bingane2022e} constructed a family of convex small $n$-gons $\geo{B}_n$, for $n=2^s$ with $s\ge 3$, such that
\[
\begin{aligned}
	\ub{L}_n - L(\geo{B}_n) &= \frac{\pi^7}{32n^6} + O\left(\frac{1}{n^8}\right),\\
	\ub{W}_n - W(\geo{B}_n) &= \frac{\pi^4}{8n^4} + O\left(\frac{1}{n^6}\right).
\end{aligned}
\]
By contrast,
\[
\begin{aligned}
	\ub{L}_n - L(\geo{R}_n) &= \frac{\pi^3}{8n^2} + O \left(\frac{1}{n^4}\right),\\
	\ub{W}_n - W(\geo{R}_n) &= \frac{3\pi^2}{8n^2} + O \left(\frac{1}{n^4}\right)
\end{aligned}
\]
for all even $n\ge 4$. Note that $W(\geo{B}_8) = W_8^*$. Some polygons~$\geo{B}_n$ are illustrated in Figure~\ref{figure:Bn}. The diameter graph of~$\geo{B}_n$ has the vertical edge as axis of symmetry and can be described by a cycle of length $n/2+1$, plus $n/2-1$ additional pendant edges, arranged so that all but two particular vertices of the cycle have a pendant edge.

\begin{figure}[H]
	\centering
	\subfloat[$(\geo{B}_8,3.121062,0.977609)$]{
		\begin{tikzpicture}[scale=5]
			\draw[dashed] (0,0) -- (0.2957,0.2043) -- (0.5000,0.5000) -- (0.4114,0.9114) -- (0,1) -- (-0.4114,0.9114) -- (-0.5000,0.5000) -- (-0.2957,0.2043) -- cycle;
			\draw[red,thick] (0,0) -- (0,1);
			\draw[blue,thick] (0,0) -- (0.4114,0.9114) -- (-0.5000,0.5000) -- (0.5000,0.5000) -- (-0.4114,0.9114) -- cycle;
			\draw (0.4114,0.9114) -- (-0.2957,0.2043);\draw (-0.4114,0.9114) -- (0.2957,0.2043);
		\end{tikzpicture}
	}
	\subfloat[$(\geo{B}_{16},3.136543,0.994996)$]{
		\begin{tikzpicture}[scale=5]
			\draw[dashed] (0,0) -- (0.1838,0.0562) -- (0.3536,0.1464) -- (0.4438,0.3162) -- (0.5000,0.5000) -- (0.4800,0.6988) -- (0.3536,0.8536) -- (0.1988,0.9800) -- (0,1) -- (-0.1988,0.9800) -- (-0.3536,0.8536) -- (-0.4800,0.6988) -- (-0.5000,0.5000) -- (-0.4438,0.3162) -- (-0.3536,0.1464) -- (-0.1838,0.0562) -- cycle;
			\draw[blue,thick] (0,0) -- (0.1988,0.9800) -- (-0.3536,0.1464) -- (0.4800,0.6988) -- (-0.5000,0.5000) -- (0.5000,0.5000) -- (-0.4800,0.6988) -- (0.3536,0.1464) -- (-0.1988,0.9800) -- cycle;
			\draw[red,thick] (0,0) -- (0,1);
			\draw (0.1988,0.9800) -- (-0.1838,0.0562);\draw (-0.1988,0.9800) -- (0.1838,0.0562);
			\draw (-0.3536,0.1464) -- (0.3536,0.8536);\draw (0.3536,0.1464) -- (-0.3536,0.8536);
			\draw (0.4800,0.6988) -- (-0.4438,0.3162);\draw (-0.4800,0.6988) -- (0.4438,0.3162);
		\end{tikzpicture}
	}
	\subfloat[$(\geo{B}_{32},3.140331,0.998784)$]{
		\begin{tikzpicture}[scale=5]
			\draw[dashed] (0,0) -- (0.0966,0.0144) -- (0.1913,0.0381) -- (0.2751,0.0883) -- (0.3536,0.1464) -- (0.4117,0.2249) -- (0.4619,0.3087) -- (0.4856,0.4034) -- (0.5000,0.5000) -- (0.4951,0.5985) -- (0.4619,0.6913) -- (0.4198,0.7805) -- (0.3536,0.8536) -- (0.2805,0.9198) -- (0.1913,0.9619) -- (0.0985,0.9951) -- (0,1) -- (-0.0985,0.9951) -- (-0.1913,0.9619) -- (-0.2805,0.9198) -- (-0.3536,0.8536) -- (-0.4198,0.7805) -- (-0.4619,0.6913) -- (-0.4951,0.5985) -- (-0.5000,0.5000) -- (-0.4856,0.4034) -- (-0.4619,0.3087) -- (-0.4117,0.2249) -- (-0.3536,0.1464) -- (-0.2751,0.0883) -- (-0.1913,0.0381) -- (-0.0966,0.0144) -- cycle;
			\draw[blue,thick] (0,0) -- (0.0985,0.9951) -- (-0.1913,0.0381) -- (0.2805,0.9198) -- (-0.3536,0.1464) -- (0.4198,0.7805) -- (-0.4619,0.3087) -- (0.4951,0.5985) -- (-0.5000,0.5000) -- (0.5000,0.5000) -- (-0.4951,0.5985) -- (0.4619,0.3087) -- (-0.4198,0.7805) -- (0.3536,0.1464) -- (-0.2805,0.9198) -- (0.1913,0.0381) -- (-0.0985,0.9951) -- cycle;
			\draw[red,thick] (0,0) -- (0,1);
			\draw (0.0985,0.9951) -- (-0.0966,0.0144);\draw (-0.0985,0.9951) -- (0.0966,0.0144);
			\draw (-0.1913,0.0381) -- (0.1913,0.9619);\draw (0.1913,0.0381) -- (-0.1913,0.9619);
			\draw (0.2805,0.9198) -- (-0.2751,0.0883);\draw (-0.2805,0.9198) -- (0.2751,0.0883);
			\draw (-0.3536,0.1464) -- (0.3536,0.8536);\draw (0.3536,0.1464) -- (-0.3536,0.8536);
			\draw (0.4198,0.7805) -- (-0.4117,0.2249);\draw (-0.4198,0.7805) -- (0.4117,0.2249);
			\draw (-0.4619,0.3087) -- (0.4619,0.6913);\draw (0.4619,0.3087) -- (-0.4619,0.6913);
			\draw (0.4951,0.5985) -- (-0.4856,0.4034);\draw (-0.4951,0.5985) -- (0.4856,0.4034);
		\end{tikzpicture}
	}
	\caption{Best prior polygons $(\geo{B}_n,L(\geo{B}_n),W(\geo{B}_n))$: (a) Octagon $\geo{B}_8$; (b) Hexadecagon $\geo{B}_{16}$; (c) Triacontadigon $\geo{B}_{32}$}
	\label{figure:Bn}
\end{figure}

\section{General construction}\label{sec:Cn}
For any $n=2^s$ where $s\ge 2$ is an integer, consider a convex small $n$-gon $\geo{P}_n$ with vertices $\geo{v}_0,\geo{v}_1,\ldots,\geo{v}_{n-1}$ having the following diameter graph:
\begin{itemize}
	\item a cycle of odd length $m$: $\geo{v}_{0} - \geo{v}_1 - \ldots-\geo{v}_{k} - \ldots - \geo{v}_{\frac{m-1}{2}} - \geo{v}_{\frac{m+1}{2}} - \ldots - \geo{v}_{m-k} - \ldots - \geo{v}_{m-1}-\geo{v}_0$,
	\item the pendant edge $\geo{v}_{0}-\geo{v}_{m}$ as axis of symmetry, and
	\item $n-m-1$ other pendant edges from some particular vertices $\geo{v}_k$ of the cycle,
\end{itemize}
as illustrated in Figure~\ref{figure:model}.

\begin{figure}[H]
	\centering
	\begin{tikzpicture}[scale=10]
		\draw[dashed] (0,0) node[below]{$\geo{v}_0(0,0)$} -- (0.1860,0.0566) node[below right]{$\geo{v}_{15}(x_{15},y_{15})$} -- (0.3576,0.1481) node[below right]{$\geo{v}_9(x_9,y_9)$} -- (0.4811,0.2983) node[right]{$\geo{v}_7(x_7,y_7)$} -- (0.5000,0.4919) node[right]{$\geo{v}_5(x_5,y_5)$} -- (0.4428,0.6810) node[right]{$\geo{v}_{13}(x_{13},y_{13})$} -- (0.3495,0.8552) node[right]{$\geo{v}_3(x_3,y_3)$} -- (0.1966,0.9805) node[above right]{$\geo{v}_1(x_1,y_1)$} -- (0,1) node[above]{$\geo{v}_{11}(0,1)$} -- (-0.1966,0.9805) node[above left]{$\geo{v}_{10}(x_{10},y_{10})$} -- (-0.3495,0.8552) node[left]{$\geo{v}_8(x_8,y_8)$} -- (-0.4428,0.6810) node[left]{$\geo{v}_{14}(x_{14},y_{14})$} -- (-0.5000,0.4919) node[left]{$\geo{v}_6(x_6,y_6)$} -- (-0.4811,0.2983) node[left]{$\geo{v}_4(x_4,y_4)$} -- (-0.3576,0.1481) node[below left]{$\geo{v}_2(x_2,y_2)$} -- (-0.1860,0.0566) node[below left]{$\geo{v}_{12}(x_{12},y_{12})$} -- cycle;
		\draw (0,0)--(0,1);
		\draw (0.5000,0.4919)--(-0.5000,0.4919);
		\draw (0,0) -- (0.1966,0.9805) -- (-0.3576,0.1481) -- (0.3495,0.8552) -- (-0.4847,0.3006) -- (0.5000,0.4919);
		\draw (0,0) -- (-0.1966,0.9805) -- (0.3576,0.1481) -- (-0.3495,0.8552) -- (0.4847,0.3006) -- (-0.5000,0.4919);
		\draw (0.1966,0.9805) -- (-0.1860,0.0566);\draw (-0.1966,0.9805) -- (0.1860,0.0566);
		\draw (-0.4847,0.3006) -- (0.4428,0.6810);\draw (0.4847,0.3006) -- (-0.4428,0.6810);
		\draw (0.0492,0.2451) arc (78.66:90.00:0.25) node[midway,above]{${\color{red}\alpha_0}$};
		\draw (0.1242,0.7425) node{${\color{red}\alpha_1}$};
		\draw (0.0795,0.7610) node{${\color{red}\alpha_1}$};
		\draw (-0.1999,0.3406) node{${\color{red}\alpha_2}$};
		\draw (0.1573,0.6972) node{${\color{red}\alpha_3}$};
		\draw (-0.2618,0.4158) node{${\color{red}\alpha_4}$};
		\draw (-0.2430,0.3704) node{${\color{red}\alpha_4}$};
		\draw (0.2524,0.4677) node{${\color{red}\alpha_5}$};
	\end{tikzpicture}
	\caption{Definition of variables $\alpha_0, \alpha_1, \ldots, \alpha_{\frac{m-1}{2}}$: Case of $n=16$  and $m=11$}
	\label{figure:model}
\end{figure}

Let $c_0 \alpha_0 := \angle \geo{v}_{m} \geo{v}_{0} \geo{v}_{1}$, and for $k = 1,2,\ldots, \frac{m-1}{2}$, let $c_k\alpha_k := \angle \geo{v}_{k-1} \geo{v}_{k} \geo{v}_{k+1}$, where $c_0$ and $c_k$ are the numbers of sides of $\geo{P}_n$ seen from the angles $\angle \geo{v}_{m} \geo{v}_{0} \geo{v}_{1}$ and $\angle \geo{v}_{k-1} \geo{v}_{k} \geo{v}_{k+1}$, respectively. Due to the symmetry of the construction, we can describe the $n$-gon $\geo{P}_n$ as a composition $(c_0,c_1,\ldots,c_{\frac{m-1}{2}})$ of $\frac{n}{2}$ into $\frac{m+1}{2}$ parts, i.e.,
\[
\sum_{k=0}^{\frac{m-1}{2}} c_k = \frac{n}{2}.
\]
For example, for the polygons $\geo{B}_n$ in Figure~\ref{figure:Bn}, we have $c_0=1$, $c_1 =\ldots = c_{n/4-1}=2$, $c_{n/4}=1$. One can remark that the sequence $(c_k)$ corresponding to~$\geo{B}_n$ consists of the digits of the palindromic number $11\prod_{j=3}^{\log_2 n}(1+10^{n/2^j})$.

Let us assume that, if $c_k > 1$, the pendant edges from the vertex $\geo{v}_k$ subdivide the angle $c_k\alpha_k$ into $c_k$~equal angles $\alpha_k$. We then have

\begin{equation}\label{eq:condition}
	\sum_{k=0}^{\frac{m-1}{2}} c_k \alpha_k = \frac{\pi}{2},
\end{equation}
and
\begin{subequations}\label{eq:LW}
	\begin{align}
		L(\geo{P}_n) &= \sum_{k=0}^{\frac{m-1}{2}} 4c_k\sin \frac{\alpha_k}{2},\label{eq:LW:L}\\
		W(\geo{P}_n) &= \min_{k=0,1,\ldots,\frac{m-1}{2}} \cos \frac{\alpha_k}{2}.
	\end{align}
\end{subequations}

We use cartesian coordinates to describe the $n$-gon $\geo{P}_n$, assuming that a vertex $\geo{v}_k$ is positioned at abscissa $x_k$ and ordinate $y_k$. Placing the vertex $\geo{v}_0$ at the origin and the vertex $\geo{v}_{m}$ at $(0,1)$ in the plane, we have
\begin{subequations}\label{eq:xy}
	\begin{align}
		x_{k} &= \sum_{j=1}^{k} (-1)^{j-1} \sin \left(\sum_{i=0}^{j-1} c_i \alpha_i\right) && = - x_{m-k} &\forall k=1,2,\ldots, \frac{m-1}{2},\\
		y_{k} &= \sum_{j=1}^{k} (-1)^{j-1} \cos \left(\sum_{i=0}^{j-1} c_i \alpha_i\right) && = y_{m-k} &\forall k=1,2,\ldots, \frac{m-1}{2}.
	\end{align}
\end{subequations}
Since the edge $\geo{v}_{\frac{m-1}{2}} - \geo{v}_{\frac{m+1}{2}}$ is horizontal and $\|\geo{v}_{\frac{m-1}{2}} - \geo{v}_{\frac{m+1}{2}}\| = 1$, we also have
\begin{equation}\label{eq:x4}
	x_{\frac{m-1}{2}} = (-1)^{\frac{m+1}{2}}/2.
\end{equation}

In~\cite{mossinghoff2006b,bingane2022e}, it is observed that the angles $\alpha_k$ that maximize $L(\geo{P}_n)$ in~\eqref{eq:LW:L} exhibit a pattern of damped oscillation, converging in an alterning manner to a mean value around $\pi/n$. To mimic this pattern, we suppose
\begin{equation}
	\alpha_k =
	\begin{cases}
		\frac{\pi}{n} + \delta_0 &\text{if $k$ is even,}\\
		\frac{\pi}{n} - \delta_1 &\text{if $k$ is odd.}
	\end{cases}
\end{equation}
Then \eqref{eq:condition} and \eqref{eq:LW} become
\begin{equation}\label{eq:conditiond}
	\sum_{\text{even $k$}} c_k \delta_0 = \sum_{\text{odd $k$}} c_k \delta_1
\end{equation}
and
\begin{subequations}\label{eq:LWd}
	\begin{align}
		L(\geo{P}_n) &= \sum_{\text{even $k$}} 4c_{k} \sin \left(\frac{\pi}{2n} + \frac{\delta_0}{2}\right) + \sum_{\text{odd $k$}} 4c_{k} \sin \left(\frac{\pi}{2n} - \frac{\delta_1}{2}\right), \label{eq:LWd:L}\\
		W(\geo{P}_n) &= \min \left\{\cos \left(\frac{\pi}{2n} + \frac{\delta_0}{2}\right), \cos \left(\frac{\pi}{2n} - \frac{\delta_1}{2}\right)\right\},
	\end{align}
\end{subequations}
respectively. Using~\eqref{eq:conditiond} to eliminate $\delta_1$, \eqref{eq:x4} becomes an equation
\begin{equation}
	\label{eq:x4d}
	x_{\frac{m-1}{2}}(\delta_0) = (-1)^{\frac{m+1}{2}}/2
\end{equation}
that relies only on $\delta_0$. For the polygons $\geo{B}_n$ in Figure~\ref{figure:Bn}, this equation becomes
\[
\frac{\sin \left(\delta_0 -\frac{\pi}{n}\right)}{\sin \frac{2\pi}{n}} = -\frac{1}{2} \Rightarrow \delta_0 = \delta_1 = \frac{\pi}{n}-\arcsin \left(\frac{1}{2}\sin \frac{2\pi}{n}\right) = \frac{\pi^3}{2n^3} + O\left(\frac{1}{n^5}\right).
\]

We can now prove Theorem~\ref{thm:Cn}.

\begin{proof}[Proof of Theorem~\ref{thm:Cn}]
	Choosing $m = 3n/4-1$ and
	\[
	c_k =
	\begin{cases}
		2 &\text{if $k = 3j-2$,}\\
		1 &\text{otherwise,}
	\end{cases}
	\]
	Equation~\eqref{eq:x4d} becomes
	\[
	\frac{\cos \frac{2\pi}{n} + \sin \frac{2\pi}{n}}{2\cos \frac{2\pi}{n}} - \frac{\sin \frac{\pi}{n} \cos\delta_0}{\cos \frac{2\pi}{n}} + \frac{\cos \frac{\pi}{n} \sin \delta_0}{\sin \frac{2\pi}{n}} = \frac{1}{2}
	\]
	and has a solution $\delta_0 = \delta_1 = \delta(n)$ satisfying
	\[
	\delta(n) = \arctan \left(\sec \frac{2\pi}{n}-1\right) - \arcsin\left(\sin \left(\arctan \left(\sec \frac{2\pi}{n}-1\right)\right)\cos \frac{\pi}{n}\right) = \frac{\pi^4}{n^4} + O\left(\frac{1}{n^6}\right).
	\]
	Let $\geo{C}_n$ denote the $n$-gon obtained by setting $\delta_0 = \delta(n)$. We have, from~\eqref{eq:LWd},
	\[
	\begin{aligned}
		\ub{L}_n - L(\geo{C}_n) &= \frac{\pi^9}{8n^8} + O\left(\frac{1}{n^{10}}\right),\\
		\ub{W}_n - W(\geo{C}_n) &= \frac{\pi^5}{4n^5} + O\left(\frac{1}{n^7}\right).
	\end{aligned}
	\]
	We illustrate $\geo{C}_n$ for some $n$ in Figure~\ref{figure:Cn}.
\end{proof}

\begin{figure}[h]
	\centering
	\subfloat[$(\geo{C}_{16},3.136548,0.995107)$]{
		\begin{tikzpicture}[scale=5]
			\draw[dashed] (0,0) -- (0.1860,0.0566) -- (0.3576,0.1481) -- (0.4811,0.2983) -- (0.5000,0.4919) -- (0.4428,0.6810) -- (0.3495,0.8552) -- (0.1966,0.9805) -- (0,1) -- (-0.1966,0.9805) -- (-0.3495,0.8552) -- (-0.4428,0.6810) -- (-0.5000,0.4919) -- (-0.4811,0.2983) -- (-0.3576,0.1481) -- (-0.1860,0.0566) -- cycle;
			\draw[red,thick] (0,0)--(0,1);
			\draw[blue,thick] (0.5000,0.4919)--(-0.5000,0.4919);
			\draw[blue,thick] (0,0) -- (0.1966,0.9805) -- (-0.3576,0.1481) -- (0.3495,0.8552) -- (-0.4847,0.3006) -- (0.5000,0.4919);
			\draw[blue,thick] (0,0) -- (-0.1966,0.9805) -- (0.3576,0.1481) -- (-0.3495,0.8552) -- (0.4847,0.3006) -- (-0.5000,0.4919);
			\draw (0.1966,0.9805) -- (-0.1860,0.0566);\draw (-0.1966,0.9805) -- (0.1860,0.0566);
			\draw (-0.4847,0.3006) -- (0.4428,0.6810);\draw (0.4847,0.3006) -- (-0.4428,0.6810);
		\end{tikzpicture}
	}
	\subfloat[$(\geo{C}_{32},3.140331,0.998793)$]{
		\begin{tikzpicture}[scale=5]
			\draw[dashed] (0,0) -- (0.0970,0.0144) -- (0.1921,0.0382) -- (0.2807,0.0801) -- (0.3534,0.1460) -- (0.4118,0.2247) -- (0.4622,0.3088) -- (0.4952,0.4011) -- (0.5000,0.4990) -- (0.4856,0.5962) -- (0.4617,0.6915) -- (0.4197,0.7803) -- (0.3538,0.8531) -- (0.2749,0.9116) -- (0.1906,0.9621) -- (0.0981,0.9952) -- (0,1) -- (-0.0981,0.9952) -- (-0.1906,0.9621) -- (-0.2749,0.9116) -- (-0.3538,0.8531) -- (-0.4197,0.7803) -- (-0.4617,0.6915) -- (-0.4856,0.5962) -- (-0.5000,0.4990) -- (-0.4952,0.4011) -- (-0.4622,0.3088) -- (-0.4118,0.2247) -- (-0.3534,0.1460) -- (-0.2807,0.0801) -- (-0.1921,0.0382) -- (-0.0970,0.0144) -- cycle;
			\draw[red,thick] (0,0)--(0,1);
			\draw[blue,thick] (0.5000,0.4990)--(-0.5000,0.4990);
			\draw[blue,thick] (0,0) -- (0.0981,0.9952) -- (-0.1921,0.0382) -- (0.1906,0.9621) -- (-0.2807,0.0801) -- (0.3538,0.8531) -- (-0.3534,0.1460) -- (0.4197,0.7803) -- (-0.4622,0.3088) -- (0.4617,0.6915) -- (-0.4952,0.4011) -- (0.5000,0.4990);
			\draw[blue,thick] (0,0) -- (-0.0981,0.9952) -- (0.1921,0.0382) -- (-0.1906,0.9621) -- (0.2807,0.0801) -- (-0.3538,0.8531) -- (0.3534,0.1460) -- (-0.4197,0.7803) -- (0.4622,0.3088) -- (-0.4617,0.6915) -- (0.4952,0.4011) -- (-0.5000,0.4990);
			\draw (0.0981,0.9952) -- (-0.0970,0.0144);\draw (-0.0981,0.9952) -- (0.0970,0.0144);
			\draw (-0.2807,0.0801) -- (0.2749,0.9116);\draw (0.2807,0.0801) -- (-0.2749,0.9116);
			\draw (0.4197,0.7803) -- (-0.4118,0.2247);\draw (-0.4197,0.7803) -- (0.4118,0.2247);
			\draw (-0.4952,0.4011) -- (0.4856,0.5962);\draw (0.4952,0.4011) -- (-0.4856,0.5962);
		\end{tikzpicture}
	}
	\subfloat[$(\geo{C}_{64},3.141277,0.999699)$]{
		\begin{tikzpicture}[scale=5]
			\draw[dashed] (0,0) -- (0.0489,0.0036) -- (0.0977,0.0096) -- (0.1455,0.0204) -- (0.1913,0.0380) -- (0.2352,0.0601) -- (0.2779,0.0843) -- (0.3180,0.1125) -- (0.3535,0.1464) -- (0.3856,0.1835) -- (0.4158,0.2222) -- (0.4420,0.2637) -- (0.4619,0.3086) -- (0.4773,0.3552) -- (0.4904,0.4025) -- (0.4988,0.4508) -- (0.5000,0.4999) -- (0.4964,0.5488) -- (0.4904,0.5976) -- (0.4796,0.6455) -- (0.4620,0.6913) -- (0.4399,0.7351) -- (0.4157,0.7778) -- (0.3874,0.8179) -- (0.3536,0.8535) -- (0.3164,0.8856) -- (0.2777,0.9158) -- (0.2362,0.9420) -- (0.1914,0.9619) -- (0.1447,0.9773) -- (0.0974,0.9904) -- (0.0491,0.9988) -- (0,1) -- (-0.0491,0.9988) -- (-0.0974,0.9904) -- (-0.1447,0.9773) -- (-0.1914,0.9619) -- (-0.2362,0.9420) -- (-0.2777,0.9158) -- (-0.3164,0.8856) -- (-0.3536,0.8535) -- (-0.3874,0.8179) -- (-0.4157,0.7778) -- (-0.4399,0.7351) -- (-0.4620,0.6913) -- (-0.4796,0.6455) -- (-0.4904,0.5976) -- (-0.4964,0.5488) -- (-0.5000,0.4999) -- (-0.4988,0.4508) -- (-0.4904,0.4025) -- (-0.4773,0.3552) -- (-0.4619,0.3086) -- (-0.4420,0.2637) -- (-0.4158,0.2222) -- (-0.3856,0.1835) -- (-0.3535,0.1464) -- (-0.3180,0.1125) -- (-0.2779,0.0843) -- (-0.2352,0.0601) -- (-0.1913,0.0380) -- (-0.1455,0.0204) -- (-0.0977,0.0096) -- (-0.0489,0.0036) -- cycle;
			\draw[red,thick] (0,0)--(0,1);
			\draw[blue,thick] (0.5000,0.5000)--(-0.5000,0.4999);
			\draw[blue,thick] (0,0) -- (0.0491,0.9988) -- (-0.0977,0.0096) -- (0.0974,0.9904) -- (-0.1455,0.0204) -- (0.1914,0.9619) -- (-0.1913,0.0380) -- (0.2362,0.9420) -- (-0.2779,0.0843) -- (0.2777,0.9158) -- (-0.3180,0.1125) -- (0.3536,0.8535) -- (-0.3535,0.1464) -- (0.3874,0.8179) -- (-0.4158,0.2222) -- (0.4157,0.7778) -- (-0.4420,0.2637) -- (0.4620,0.6913) -- (-0.4619,0.3086) -- (0.4796,0.6455) -- (-0.4904,0.4025) -- (0.4904,0.5976) -- (-0.4988,0.4508) -- (0.5000,0.4999);
			\draw[blue,thick] (0,0) -- (-0.0491,0.9988) -- (0.0977,0.0096) -- (-0.0974,0.9904) -- (0.1455,0.0204) -- (-0.1914,0.9619) -- (0.1913,0.0380) -- (-0.2362,0.9420) -- (0.2779,0.0843) -- (-0.2777,0.9158) -- (0.3180,0.1125) -- (-0.3536,0.8535) -- (0.3535,0.1464) -- (-0.3874,0.8179) -- (0.4158,0.2222) -- (-0.4157,0.7778) -- (0.4420,0.2637) -- (-0.4620,0.6913) -- (0.4619,0.3086) -- (-0.4796,0.6455) -- (0.4904,0.4025) -- (-0.4904,0.5976) -- (0.4988,0.4508) -- (-0.5000,0.4999);
			\draw (0.0491,0.9988) -- (-0.0489,0.0036);\draw (-0.0491,0.9988) -- (0.0489,0.0036);
			\draw (-0.1455,0.0204) -- (0.1447,0.9773);\draw (0.1455,0.0204) -- (-0.1447,0.9773);
			\draw (0.2362,0.9420) -- (-0.2352,0.0601);\draw (-0.2362,0.9420) -- (0.2352,0.0601);
			\draw (-0.3180,0.1125) -- (0.3164,0.8856);\draw (0.3180,0.1125) -- (-0.3164,0.8856);
			\draw (0.3874,0.8179) -- (-0.3856,0.1835);\draw (-0.3874,0.8179) -- (0.3856,0.1835);
			\draw (-0.4420,0.2637) -- (0.4399,0.7351);\draw (0.4420,0.2637) -- (-0.4399,0.7351);
			\draw (0.4796,0.6455) -- (-0.4773,0.3552);\draw (-0.4796,0.6455) -- (0.4773,0.3552);
			\draw (-0.4988,0.4508) -- (0.4964,0.5488);\draw (0.4988,0.4508) -- (-0.4964,0.5488);
		\end{tikzpicture}
	}
	\caption{Polygons $(\geo{C}_n,L(\geo{C}_n),W(\geo{C}_n))$ defined in Theorem~\ref{thm:Cn}: (a) Hexadecagon~$\geo{C}_{16}$; (b) Triacontadigon~$\geo{C}_{32}$; (c) Hexacontatetragon~$\geo{C}_{64}$}
	\label{figure:Cn}
\end{figure}

For $n=2^s \ge 16$, the $n$-gon $\geo{C}_n$ provides tighter lower bounds on the maximal perimeter $L_n^*$ and the maximal width $W_n^*$ compared to $\geo{B}_n$. For a numerical comparison, the perimeters and the widths of $\geo{B}_n$ and $\geo{C}_n$ are given in Table~\ref{table:L(Pn)}.

For $n=2^s \ge 32$, the lower bounds provided by $\geo{C}_n$ can still be improved. Using Lemma~\ref{thm:ps}, a better family of convex small $n$-gons $\geo{J}_n$, for $n=2^s$ with $s\ge 2$, proposed by Xue, Lian, Wang, and Zhang~\cite{xue2021} is defined in Theorem~\ref{thm:Jn}. For each $n$, the diameter graph of $\geo{J}_n$ has a cycle of length $\lfloor 2n/3 \rceil$. The family of $n$-gons $\geo{J}_n$ includes the $4$-gon of maximal perimeter and maximal width $\geo{T}_4$, the $8$-gon of maximal width $\geo{B}_8$ and the $16$-gon $\geo{C}_{16}$.

\begin{lemma}
	\label{thm:ps}
	For a positive integer $s \ge 2$, let $p_s := \prod_{j=1}^{s-1}(1+10^{\lfloor 2^j/3\rceil})$.
	\begin{enumerate}
		\item $p_s$ has exactly $\lambda_s := \frac{2^{s+1}+3+(-1)^s}{6}$ digits.
		\item Let $c_0,c_1,\ldots,c_{\lambda_s-1}$ be the digits of $p_s$, i.e., $p_s =\sum_{k=0}^{\lambda_s-1} 10^kc_k$. Then
		\begin{enumerate}
			\item $c_0,c_1,\ldots,c_{\lambda_s-1} \in \{1,2\}$ with $c_0 = c_{\lambda_s-1} = 1$,
			\item $\sum_{k=0}^{\lambda_s-1} c_k = 2^{s-1}$ and $\sum_{k=0}^{\lambda_s-1} (-1)^kc_k = 0$,
			\item $p_s$ is a palindromic number.
		\end{enumerate}
	\end{enumerate}
\end{lemma}
\begin{proof}
	The proof can be done by induction and it is left to the reader.
\end{proof}

\begin{theorem}[Xue, Lian, Wang, and Zhang~\cite{xue2021}]
	\label{thm:Jn}
	Suppose $n=2^s$ with integer $s\ge 2$. Then there exists a convex small $n$-gon $\geo{J}_n$ such that
	\[
	\begin{aligned}
		\ub{L}_n - L(\geo{J}_n) &\sim \frac{\sigma^2\pi^{2s+5}}{32n^{s+5}},\\
		\ub{W}_n - W(\geo{J}_n) &\sim \frac{\sigma\pi^{s+3}}{8n^{\frac{s+7}{2}}},
	\end{aligned}
	\]
	where $\sigma := \prod_{j=2}^\infty \frac{\sin (\pi/2^j)}{\pi/2^j} = 0.869861\ldots$.
\end{theorem}
\begin{proof}
	Let $m = \frac{2n+(-1)^s}{3}$, which is a Jacobsthal number (sequence \href{https://oeis.org/A001045}{A001045} in the \href{https://oeis.org/}{OEIS}) and let $c_0,c_1,\ldots,c_{\frac{m-1}{2}}$ be the digits of the palindromic number $p_s = \prod_{k=1}^{s-1}(1+10^{\lfloor 2^k/3\rceil})$ defined in Lemma~\ref{thm:ps}. We have $\sum_{k=0}^{\frac{m-1}{2}}c_k = n/2$ and $\sum_{k=0}^{\frac{m-1}{2}}(-1)^k c_k = 0$.
	
	After some simplifications, Equation~\eqref{eq:x4d} can be rewritten as
	\[
	a_s \left(\cos \delta_0 - \cos \frac{\pi}{n}\right) + \frac{\sin \delta_0}{2\sin \frac{\pi}{n}} = 0
	\]
	with $a_s = (-1)^{\lceil s/2 \rceil-1}2^{s-2}\prod_{j=2}^{s} \sin \frac{\pi}{2^j}$. This equation has a solution $\delta_0 = \delta_1 = \delta (n)$ satisfying
	\[
	\delta(n) = -\arctan \left(2a_s\sin \frac{\pi}{n}\right)+\arcsin\left(\sin\left(\arctan \left(2a_s\sin \frac{\pi}{n}\right)\right)\cos\frac{\pi}{n}\right) \sim -a_s\frac{\pi^3}{n^3}.
	\]
	
	Let $\geo{J}_n$ denote the $n$-gon obtained by setting $\delta_0 = \delta(n)$. We have, from~\eqref{eq:LWd},
	\[
	\begin{aligned}
		L(\geo{J}_n) &= 2n\sin \frac{\pi}{2n} \cos \frac{\delta(n)}{2},\\
		W(\geo{J}_n) &= \cos \left(\frac{\pi}{2n} + \frac{|\delta(n)|}{2}\right).
	\end{aligned}
	\]
	If $\sigma := \prod_{j=2}^\infty \frac{\sin (\pi/2^j)}{\pi/2^j} = 0.869861\ldots$, then
	\[
	a_s \sim (-1)^{\lceil s/2 \rceil-1}2^{s-2}\sigma \prod_{j=2}^{s} \frac{\pi}{2^j} = (-1)^{\lceil s/2 \rceil-1}\frac{\sigma \pi^{s-1}}{2n^{\frac{s-1}{2}}} \Rightarrow \delta(n)\sim (-1)^{\lceil s/2 \rceil}\frac{\sigma \pi^{s+2}}{2n^{\frac{s+5}{2}}}
	\]
	and
	\[
	\begin{aligned}
		\ub{L}_n - L(\geo{J}_n) &\sim \frac{\pi\delta^2(n)}{8} \sim \frac{\sigma^2\pi^{2s+5}}{32n^{s+5}},\\
		\ub{W}_n - W(\geo{J}_n) &\sim \frac{\pi|\delta(n)|}{4n} \sim \frac{\sigma\pi^{s+3}}{8n^{\frac{s+7}{2}}}.
	\end{aligned}
	\]
	We illustrate $\geo{J}_n$ for some $n$ in Figure~\ref{figure:Jn}.
\end{proof}

\begin{figure}[h]
	\centering
	\subfloat[$(\geo{J}_4,3.035276,0.866025)$]{
		\begin{tikzpicture}[scale=5]
			\draw[dashed] (0.5000,0.8660) -- (0,1) -- (-0.5000,0.8660);
			\draw[blue,thick] (0,0) -- (0.5000,0.8660) -- (-0.5000,0.8660) -- cycle;
			\draw[red,thick] (0,0) -- (0,1);
		\end{tikzpicture}
	}
	\subfloat[$(\geo{J}_8,3.121062,0.977609)$]{
		\begin{tikzpicture}[scale=5]
			\draw[dashed] (0,0) -- (0.2957,0.2043) -- (0.5000,0.5000) -- (0.4114,0.9114) -- (0,1) -- (-0.4114,0.9114) -- (-0.5000,0.5000) -- (-0.2957,0.2043) -- cycle;
			\draw[blue,thick] (0,0) -- (0.4114,0.9114) -- (-0.5000,0.5000) -- (0.5000,0.5000) -- (-0.4114,0.9114) -- cycle;
			\draw[red,thick] (0,0) -- (0,1);
			\draw (0.4114,0.9114) -- (-0.2957,0.2043);\draw (-0.4114,0.9114) -- (0.2957,0.2043);
		\end{tikzpicture}
	}
	\subfloat[$(\geo{J}_{16},3.136548,0.995107)$]{
		\begin{tikzpicture}[scale=5]
			\draw[dashed] (0,0) -- (0.1860,0.0566) -- (0.3576,0.1481) -- (0.4811,0.2983) -- (0.5000,0.4919) -- (0.4428,0.6810) -- (0.3495,0.8552) -- (0.1966,0.9805) -- (0,1) -- (-0.1966,0.9805) -- (-0.3495,0.8552) -- (-0.4428,0.6810) -- (-0.5000,0.4919) -- (-0.4811,0.2983) -- (-0.3576,0.1481) -- (-0.1860,0.0566) -- cycle;
			\draw[red,thick] (0,0)--(0,1);
			\draw[blue,thick] (0.5000,0.4919)--(-0.5000,0.4919);
			\draw[blue,thick] (0,0) -- (0.1966,0.9805) -- (-0.3576,0.1481) -- (0.3495,0.8552) -- (-0.4847,0.3006) -- (0.5000,0.4919);
			\draw[blue,thick] (0,0) -- (-0.1966,0.9805) -- (0.3576,0.1481) -- (-0.3495,0.8552) -- (0.4847,0.3006) -- (-0.5000,0.4919);
			\draw (0.1966,0.9805) -- (-0.1860,0.0566);\draw (-0.1966,0.9805) -- (0.1860,0.0566);
			\draw (-0.4847,0.3006) -- (0.4428,0.6810);\draw (0.4847,0.3006) -- (-0.4428,0.6810);
		\end{tikzpicture}
	}\\
	\subfloat[$(\geo{J}_{32},3.140331,0.998794)$]{
		\begin{tikzpicture}[scale=7]
			\draw[dashed] (0,0) -- (0.0971,0.0144) -- (0.1923,0.0383) -- (0.2811,0.0802) -- (0.3538,0.1462) -- (0.4198,0.2189) -- (0.4617,0.3077) -- (0.4856,0.4029) -- (0.5000,0.5000) -- (0.4952,0.5980) -- (0.4621,0.6903) -- (0.4117,0.7745) -- (0.3533,0.8533) -- (0.2745,0.9117) -- (0.1903,0.9621) -- (0.0980,0.9952) -- (0,1) -- (-0.0980,0.9952) -- (-0.1903,0.9621) -- (-0.2745,0.9117) -- (-0.3533,0.8533) -- (-0.4117,0.7745) -- (-0.4621,0.6903) -- (-0.4952,0.5980) -- (-0.5000,0.5000) -- (-0.4856,0.4029) -- (-0.4617,0.3077) -- (-0.4198,0.2189) -- (-0.3538,0.1462) -- (-0.2811,0.0802) -- (-0.1923,0.0383) -- (-0.0971,0.0144) -- cycle;
			\draw[red,thick] (0,0)--(0,1);
			\draw[blue,thick] (0.5000,0.4990)--(-0.5000,0.4990);
			\draw[blue,thick] (0,0) -- (0.0980,0.9952) -- (-0.1923,0.0383) -- (0.1903,0.9621) -- (-0.2811,0.0802) -- (0.3533,0.8533) -- (-0.4198,0.2189) -- (0.4621,0.6903) -- (-0.4617,0.3077) -- (0.4952,0.5980) -- (-0.5000,0.5000);
			\draw[blue,thick] (0,0) -- (-0.0980,0.9952) -- (0.1923,0.0383) -- (-0.1903,0.9621) -- (0.2811,0.0802) -- (-0.3533,0.8533) -- (0.4198,0.2189) -- (-0.4621,0.6903) -- (0.4617,0.3077) -- (-0.4952,0.5980) -- (0.5000,0.5000);
			\draw (0.0980,0.9952) -- (-0.0971,0.0144);\draw (-0.0980,0.9952) -- (0.0971,0.0144);
			\draw (-0.2811,0.0802) -- (0.2745,0.9117);\draw (0.2811,0.0802) -- (-0.2745,0.9117);
			\draw (0.3533,0.8533) -- (-0.3538,0.1462);\draw (-0.3533,0.8533) -- (0.3538,0.1462);
			\draw (-0.4198,0.2189) -- (0.4117,0.7745);\draw (0.4198,0.2189) -- (-0.4117,0.7745);
			\draw (0.4952,0.5980) -- (-0.4856,0.4029);\draw (-0.4952,0.5980) -- (0.4856,0.4029);
		\end{tikzpicture}
	}
	\subfloat[$(\geo{J}_{64},3.141277,0.999699)$]{
		\begin{tikzpicture}[scale=7]
			\draw[dashed] (0,0) -- (0.0489,0.0036) -- (0.0977,0.0096) -- (0.1456,0.0204) -- (0.1913,0.0380) -- (0.2362,0.0579) -- (0.2777,0.0842) -- (0.3164,0.1144) -- (0.3535,0.1464) -- (0.3874,0.1820) -- (0.4157,0.2221) -- (0.4398,0.2648) -- (0.4619,0.3087) -- (0.4773,0.3553) -- (0.4904,0.4026) -- (0.4988,0.4509) -- (0.5000,0.5000) -- (0.4964,0.5490) -- (0.4904,0.5977) -- (0.4796,0.6456) -- (0.4620,0.6914) -- (0.4421,0.7362) -- (0.4158,0.7777) -- (0.3856,0.8164) -- (0.3536,0.8536) -- (0.3180,0.8874) -- (0.2779,0.9157) -- (0.2352,0.9399) -- (0.1913,0.9619) -- (0.1447,0.9773) -- (0.0974,0.9904) -- (0.0491,0.9988) -- (0,1) -- (-0.0491,0.9988) -- (-0.0974,0.9904) -- (-0.1447,0.9773) -- (-0.1913,0.9619) -- (-0.2352,0.9399) -- (-0.2779,0.9157) -- (-0.3180,0.8874) -- (-0.3536,0.8536) -- (-0.3856,0.8164) -- (-0.4158,0.7777) -- (-0.4421,0.7362) -- (-0.4620,0.6914) -- (-0.4796,0.6456) -- (-0.4904,0.5977) -- (-0.4964,0.5490) -- (-0.5000,0.5000) -- (-0.4988,0.4509) -- (-0.4904,0.4026) -- (-0.4773,0.3553) -- (-0.4619,0.3087) -- (-0.4398,0.2648) -- (-0.4157,0.2221) -- (-0.3874,0.1820) -- (-0.3535,0.1464) -- (-0.3164,0.1144) -- (-0.2777,0.0842) -- (-0.2362,0.0579) -- (-0.1913,0.0380) -- (-0.1456,0.0204) -- (-0.0977,0.0096) -- (-0.0489,0.0036) -- cycle;
			\draw[red,thick] (0,0)--(0,1);
			\draw[blue,thick] (0.5000,0.5000)--(-0.5000,0.5000);
			\draw[blue,thick] (0,0) -- (0.0491,0.9988) -- (-0.0977,0.0096) -- (0.0974,0.9904) -- (-0.1456,0.0204) -- (0.1913,0.9619) -- (-0.2362,0.0579) -- (0.2779,0.9157) -- (-0.2777,0.0842) -- (0.3180,0.8874) -- (-0.3535,0.1464) -- (0.3536,0.8536) -- (-0.3874,0.1820) -- (0.4158,0.7777) -- (-0.4157,0.2221) -- (0.4421,0.7362) -- (-0.4619,0.3087) -- (0.4796,0.6456) -- (-0.4904,0.4026) -- (0.4904,0.5977) -- (-0.4988,0.4509) -- (0.5000,0.5000);
			\draw[blue,thick] (0,0) -- (-0.0491,0.9988) -- (0.0977,0.0096) -- (-0.0974,0.9904) -- (0.1456,0.0204) -- (-0.1913,0.9619) -- (0.2362,0.0579) -- (-0.2779,0.9157) -- (0.2777,0.0842) -- (-0.3180,0.8874) -- (0.3535,0.1464) -- (-0.3536,0.8536) -- (0.3874,0.1820) -- (-0.4158,0.7777) -- (0.4157,0.2221) -- (-0.4421,0.7362) -- (0.4619,0.3087) -- (-0.4796,0.6456) -- (0.4904,0.4026) -- (-0.4904,0.5977) -- (0.4988,0.4509) -- (-0.5000,0.5000);
			\draw (0.0491,0.9988) -- (-0.0489,0.0036);\draw (-0.0491,0.9988) -- (0.0489,0.0036);
			\draw (-0.1456,0.0204) -- (0.1447,0.9773);\draw (0.1456,0.0204) -- (-0.1447,0.9773);
			\draw (0.1913,0.9619) -- (-0.1913,0.0380);\draw (-0.1913,0.9619) -- (0.1913,0.0380);
			\draw (-0.2362,0.0579) -- (0.2352,0.9399);\draw (0.2362,0.0579) -- (-0.2352,0.9399);
			\draw (-0.3180,0.8874) -- (0.3164,0.1144);\draw (0.3180,0.8874) -- (-0.3164,0.1144);
			\draw (0.3874,0.1820) -- (-0.3856,0.8164);\draw (-0.3874,0.1820) -- (0.3856,0.8164);
			\draw (-0.4421,0.7362) -- (0.4398,0.2648);\draw (0.4421,0.7362) -- (-0.4398,0.2648);
			\draw (-0.4619,0.3087) -- (0.4620,0.6914);\draw (0.4619,0.3087) -- (-0.4620,0.6914);
			\draw (-0.4796,0.6456) -- (0.4773,0.3553);\draw (0.4796,0.6456) -- (-0.4773,0.3553);
			\draw (-0.4988,0.4509) -- (0.4964,0.5490);\draw (0.4988,0.4509) -- (-0.4964,0.5490);
		\end{tikzpicture}
	}
	\caption{Polygons $(\geo{J}_n,L(\geo{J}_n),W(\geo{J}_n))$ defined in Theorem~\ref{thm:Jn}: (a) Tetragon $\geo{J}_4$; (b) Octagon $\geo{J}_8$; (c) Hexadecagon $\geo{J}_{16}$; (d) Triacontadigon $\geo{J}_{32}$; (e) Hexacontatetragon $\geo{J}_{64}$}
	\label{figure:Jn}
\end{figure}

Table~\ref{table:Pn} summarizes the characteristics of $\geo{B}_n$, $\geo{C}_n$, $\geo{J}_n$ and some other families of convex small $n$-gons: $\geo{T}_n$, $\geo{Q}_n$, $\geo{G}_n$. The $n$-gons $\geo{T}_n$, $\geo{Q}_n$ and $\geo{G}_n$ are illustrated in Figure~\ref{figure:Pn}. The family of $\geo{T}_n$ represents Tamvakis polygons~\cite{tamvakis1987}. For each family of $n$-gons $\geo{P}_n$, we provide the length $m$ of the cycle of its diameter graph, the corresponding composition $(c_k)$ of $n/2$, the angle $\delta_0$ obtained by solving~\eqref{eq:x4d}, and the gaps $\ub{L}_n-L(\geo{P}_n)$ and $\ub{W}_n-W(\geo{P}_n)$. All polygons presented in this work and in~\cite{bingane2022a,bingane2022b,bingane2022c,bingane2023a} were implemented as a package: OPTIGON~\cite{optigon}, which is freely available on GitHub. In OPTIGON, we provide Julia and MATLAB functions that give the coordinates of the vertices. One can also find an algorithm developed in~\cite{bingane2022d} to find an estimate of the maximal area of a small $n$-gon when $n \ge 6$ is even.

Table~\ref{table:L(Pn)} shows the perimeters $L(\geo{P}_n)$ of the polygons defined in Table~\ref{table:Pn}, along with the upper bounds $\ub{L}_n$ and the perimeters of the regular polygons $\geo{R}_n$. The widths $W(\geo{P}_n)$ are given in Table~\ref{table:W(Pn)} with $\ub{W}_n$ and $W(\geo{R}_n)$. In all tables and figures in this paper, we rounded each numerical value at the last displayed digit. As suggested by Theorem~\ref{thm:Cn}, $\geo{C}_n$ provides tighter lower bounds on~ the maximal perimeter~$L_n^*$ and the maximal width~$W_n^*$ compared to the best prior convex small $n$-gon~$\geo{B}_n$ when $n\ge 16$, and according to Theorem~\ref{thm:Jn}, $L(\geo{J}_n) > L(\geo{C}_n)$ and $W(\geo{J}_n) > W(\geo{C}_n)$ when $n \ge 32$. For instance, we note that
\[
\begin{aligned}
	L_{128}^* - L(\geo{J}_{128}) &< \ub{L}_{128} - L(\geo{J}_{128}) &&= 3.41\ldots \times 10^{-18}\\
	&< \ub{L}_{128} - L(\geo{C}_{128}) &&= 5.19\ldots \times 10^{-14}\\
	&< \ub{L}_{128} - L(\geo{B}_{128}) &&= 2.14\ldots \times 10^{-11},\\
	W_{128}^* - W(\geo{J}_{128}) &< \ub{W}_{128} - W(\geo{J}_{128}) &&= 1.80\ldots \times 10^{-11}\\
	&< \ub{W}_{128} - W(\geo{C}_{128}) &&= 2.22\ldots \times 10^{-9}\\
	&< \ub{W}_{128} - W(\geo{B}_{128}) &&= 4.53\ldots \times 10^{-8}.
\end{aligned}
\]

\begin{table}[H]
	\footnotesize
	\centering
	\caption{Some families of polygons}
	\label{table:Pn}
	\resizebox{\linewidth}{!}{
		\begin{tabular}{@{}lllllll@{}}
			\toprule
			$n$ & $\geo{P}_n$ & $m$ & $(c_k)$ & $\delta_0$ & $\ub{L}_n-L(\geo{P}_n)$ & $\ub{W}_n-W(\geo{P}_n)$ \\
			\midrule
			$2^s\ge 4$	& $\geo{T}_n$~\cite{tamvakis1987} &	$3$	&	$c_k =
			\begin{cases} \lfloor n/6\rceil &\text{if $k = 0$,}\\\lfloor n/3\rceil &\text{if $k=1$}\end{cases}$	&	$\sim (-1)^{s-1}\frac{2\pi}{n^2}$	&	$\sim \frac{\pi^3}{4n^4}$	&	$\sim
			\begin{cases} \frac{\pi^2}{4n^3} &\text{if $s$ is even,}\\\frac{\pi^2}{2n^3} &\text{if $s$ is odd}\end{cases}$	\\
			$2^s\ge 8$	& $\geo{B}_n$~\cite{bingane2022e} &	$n/2+1$	&	$c_k =
			\begin{cases} 1 &\text{if $k = 0,n/4$,}\\2 &\text{otherwise}\end{cases}$	&	$\sim \frac{\pi^3}{2n^3}$	&	$\sim \frac{\pi^7}{32n^6}$	&	$\sim \frac{\pi^4}{8n^4}$	\\
			$2^s\ge 4$	& $\geo{J}_n$~\cite{xue2021} &	$\lfloor 2n/3 \rceil$	&	$c_k$ are digits of $p_s$\tablefootnote{$p_s =\prod_{j=1}^{s-1}(1+10^{\lfloor 2^j/3\rceil}) = \sum_{k=0}^{\frac{m-1}{2}} 10^kc_k$ (Lemma~\ref{thm:ps})}	&	$\sim (-1)^{\lceil s/2 \rceil}\frac{\sigma \pi^{s+2}}{2n^{\frac{s+5}{2}}}$\tablefootnote{$\sigma = \prod_{j=2}^\infty \frac{\sin (\pi/2^j)}{\pi/2^j} = 0.869861\ldots$ (Theorem~\ref{thm:Jn})}	&	$\sim \frac{\sigma^2\pi^{2s+5}}{32n^{s+5}}$	&	$\sim \frac{\sigma\pi^{s+3}}{8n^{\frac{s+7}{2}}}$	\\			$2^s\ge 16$	& $\geo{C}_n$ &	$3n/4-1$	&	$c_k =
			\begin{cases} 2 &\text{if $k = 3j-2$,}\\1 &\text{otherwise}\end{cases}$	&	$\sim \frac{\pi^4}{n^4}$	&	$\sim \frac{\pi^9}{8n^8}$	&	$\sim \frac{\pi^5}{4n^5}$	\\
			$2^s\ge 8$	& $\geo{G}_n$ &	$n-3$	&	$c_k =
			\begin{cases} 2 &\text{if $k = 2\lfloor n/6\rceil-1$,}\\1 &\text{otherwise}\end{cases}$	&	$\sim (-1)^{s-1}\frac{\pi^3}{\sqrt{3}n^3}$	&	$\sim \frac{\pi^7}{24n^6}$	&	$\sim \frac{\pi^4}{4\sqrt{3}n^4}$	\\
			$2^s\ge 4$	& $\geo{Q}_n$~\cite{bingane2022e} &	$n-1$	&	$c_k=1$ for all $k$	&	$\sim -\frac{\pi^2}{2n^2}$	&	$\sim \frac{\pi^5}{32n^4}$	&	$\sim \frac{\pi^3}{8n^3}$	\\
			\bottomrule
		\end{tabular}
	}
\end{table}

\begin{table}[H]
	\footnotesize
	\centering
	\caption{Perimeters of polygons defined in Table~\ref{table:Pn}}
	\label{table:L(Pn)}
	\resizebox{\linewidth}{!}{
		\begin{tabular}{@{}rllllllll@{}}
			\toprule
			$n$ & $L(\geo{R}_n)$ & $L(\geo{T}_n)$ & $L(\geo{Q}_n)$ & $L(\geo{B}_n)$ & $L(\geo{G}_n)$ & $L(\geo{C}_n)$ & $L(\geo{J}_n)$ & $\ub{L}_n$ \\
			\midrule
			4	&	2.83	&	3.0353	&	3.0353	&	--	&	--	&	--	&	3.0353	&	3.0615	\\
			8	&	3.06	&	3.1191	&	3.1193	&	3.121062	&	3.121062	&	--	&	3.121062	&	3.121445	\\
			16	&	3.1214	&	3.136438	&	3.136406	&	3.13654277	&	3.13654053	&	3.13654751	&	3.13654751	&	3.13654849	\\
			32	&	3.1365	&	3.140323	&	3.140322	&	3.14033107	&	3.14033104	&	3.1403311535	&	3.140331156355	&	3.140331156955	\\
			64	&	3.1403	&	3.14127680	&	3.14127668	&	3.1412772496	&	3.1412772491	&	3.141277250919	&	3.14127725093268	&	3.14127725093277	\\
			128	&	3.141277	&	3.14151377	&	3.14151377	&	3.141513801123	&	3.141513801116	&	3.14151380114425	&	3.14151380114430107292	&	3.14151380114430107632	\\
			\bottomrule
		\end{tabular}
	}
\end{table}

\begin{table}[H]
	\footnotesize
	\centering
	\caption{Widths of polygons defined in Table~\ref{table:Pn}}
	\label{table:W(Pn)}
	\resizebox{\linewidth}{!}{
		\begin{tabular}{@{}rllllllll@{}}
			\toprule
			$n$ & $W(\geo{R}_n)$ & $W(\geo{T}_n)$ & $W(\geo{Q}_n)$ & $W(\geo{B}_n)$ & $W(\geo{G}_n)$ & $W(\geo{C}_n)$ & $W(\geo{J}_n)$ & $\ub{W}_n$ \\
			\midrule
			4	&	0.71	&	0.87	&	0.87	&	--	&	--	&	--	&	0.87	&	0.92	\\
			8	&	0.9239	&	0.9659	&	0.9730	&	0.9776	&	0.9776	&	--	&	0.9776	&	0.9808	\\
			16	&	0.9808	&	0.9945	&	0.9942	&	0.994996	&	0.994961	&	0.995107	&	0.995107	&	0.995185	\\
			32	&	0.9952	&	0.998630	&	0.998674	&	0.998784	&	0.998782	&	0.99879314	&	0.99879450	&	0.99879546	\\
			64	&	0.9988	&	0.999689	&	0.999684	&	0.99969809	&	0.99969797	&	0.99969875	&	0.9996988128	&	0.9996988187	\\
			128	&	0.999699	&	0.99992229	&	0.99992284	&	0.99992466	&	0.99992465	&	0.9999246996	&	0.999924701821	&	0.999924701839	\\
			\bottomrule
		\end{tabular}
	}
\end{table}

\begin{figure}
	\centering
	\subfloat[$(\geo{T}_4,3.0353,0.8660)$]{
		\begin{tikzpicture}[scale=4]
			\draw[dashed] (0.5000,0.8660) -- (0,1) -- (-0.5000,0.8660);
			\draw[red,thick] (0,0) -- (0,1);
			\draw[blue,thick] (0,0) -- (0.5000,0.8660) -- (-0.5000,0.8660) -- cycle;
		\end{tikzpicture}
	}
	\subfloat[$(\geo{T}_8,3.1191,0.9659)$]{
		\begin{tikzpicture}[scale=4]
			\draw[dashed] (0,0) -- (0.2660,0.2232) -- (0.4397,0.5240) -- (0.5000,0.8660) -- (0,1) -- (-0.5000,0.8660) -- (-0.4397,0.5240) -- (-0.2660,0.2232) -- cycle;
			\draw[blue,thick] (0,0) -- (0.5000,0.8660) -- (-0.5000,0.8660) -- cycle;
			\draw[red,thick] (0,0) -- (0,1);
			\draw (0.5000,0.8660) -- (-0.4397,0.5240);\draw (0.5000,0.8660) -- (-0.2660,0.2232);
			\draw (-0.5000,0.8660) -- (0.4397,0.5240);\draw (-0.5000,0.8660) -- (0.2660,0.2232);
		\end{tikzpicture}
	}
	\subfloat[$(\geo{T}_{16},3.1364,0.9945)$]{
		\begin{tikzpicture}[scale=4]
			\draw[dashed] (0,0) -- (0.1691,0.1229) -- (0.3090,0.2782) -- (0.4135,0.4593) -- (0.4781,0.6581) -- (0.5000,0.8660) -- (0.3420,0.9397) -- (0.1736,0.9848) -- (0,1) -- (-0.1736,0.9848) -- (-0.3420,0.9397) -- (-0.5000,0.8660) -- (-0.4781,0.6581) -- (-0.4135,0.4593) -- (-0.3090,0.2782) -- (-0.1691,0.1229) -- cycle;
			\draw[blue,thick] (0,0) -- (0.5000,0.8660) -- (-0.5000,0.8660) -- cycle;
			\draw[red,thick] (0,0) -- (0,1);
			\draw (0,0) -- (0.3420,0.9397);\draw (0,0) -- (-0.3420,0.9397);
			\draw (0,0) -- (0.1736,0.9848);\draw (0,0) -- (-0.1736,0.9848);
			\draw (0.5000,0.8660) -- (-0.4781,0.6581);\draw (-0.5000,0.8660) -- (0.4781,0.6581);
			\draw (0.5000,0.8660) -- (-0.4135,0.4593);\draw (-0.5000,0.8660) -- (0.4135,0.4593);
			\draw (0.5000,0.8660) -- (-0.3090,0.2782);\draw (-0.5000,0.8660) -- (0.3090,0.2782);
			\draw (0.5000,0.8660) -- (-0.1691,0.1229);\draw (-0.5000,0.8660) -- (0.1691,0.1229);
		\end{tikzpicture}
	}
	\subfloat[$(\geo{T}_{32},3.1403,0.9986)$]{
		\begin{tikzpicture}[scale=4]
			\draw[dashed] (0,0) -- (0.0801,0.0514) -- (0.1549,0.1103) -- (0.2237,0.1759) -- (0.2861,0.2479) -- (0.3413,0.3254) -- (0.3888,0.4078) -- (0.4284,0.4944) -- (0.4595,0.5843) -- (0.4819,0.6768) -- (0.4955,0.7710) -- (0.5000,0.8660) -- (0.4067,0.9135) -- (0.3090,0.9511) -- (0.2079,0.9781) -- (0.1045,0.9945) -- (0,1) -- (-0.1045,0.9945) -- (-0.2079,0.9781) -- (-0.3090,0.9511) -- (-0.4067,0.9135) -- (-0.5000,0.8660) -- (-0.4955,0.7710) -- (-0.4819,0.6768) -- (-0.4595,0.5843) -- (-0.4284,0.4944) -- (-0.3888,0.4078) -- (-0.3413,0.3254) -- (-0.2861,0.2479) -- (-0.2237,0.1759) -- (-0.1549,0.1103) -- (-0.0801,0.0514) -- cycle;
			\draw[blue,thick] (0,0) -- (0.5000,0.8660) -- (-0.5000,0.8660) -- cycle;
			\draw[red,thick] (0,0) -- (0,1);
			\draw (0,0) -- (0.1045,0.9945);\draw (0,0) -- (-0.1045,0.9945);
			\draw (0,0) -- (0.2079,0.9781);\draw (0,0) -- (-0.2079,0.9781);
			\draw (0,0) -- (0.3090,0.9511);\draw (0,0) -- (-0.3090,0.9511);
			\draw (0,0) -- (0.4067,0.9135);\draw (0,0) -- (-0.4067,0.9135);
			\draw (0.5000,0.8660) -- (-0.4955,0.7710);\draw (-0.5000,0.8660) -- (0.4955,0.7710);
			\draw (0.5000,0.8660) -- (-0.4819,0.6768);\draw (-0.5000,0.8660) -- (0.4819,0.6768);
			\draw (0.5000,0.8660) -- (-0.4595,0.5843);\draw (-0.5000,0.8660) -- (0.4595,0.5843);
			\draw (0.5000,0.8660) -- (-0.4284,0.4944);\draw (-0.5000,0.8660) -- (0.4284,0.4944);
			\draw (0.5000,0.8660) -- (-0.3888,0.4078);\draw (-0.5000,0.8660) -- (0.3888,0.4078);
			\draw (0.5000,0.8660) -- (-0.3413,0.3254);\draw (-0.5000,0.8660) -- (0.3413,0.3254);
			\draw (0.5000,0.8660) -- (-0.2861,0.2479);\draw (-0.5000,0.8660) -- (0.2861,0.2479);
			\draw (0.5000,0.8660) -- (-0.2237,0.1759);\draw (-0.5000,0.8660) -- (0.2237,0.1759);
			\draw (0.5000,0.8660) -- (-0.1549,0.1103);\draw (-0.5000,0.8660) -- (0.1549,0.1103);
			\draw (0.5000,0.8660) -- (-0.0801,0.0514);\draw (-0.5000,0.8660) -- (0.0801,0.0514);
		\end{tikzpicture}
	}\\
	\subfloat[$(\geo{Q}_4,3.0353,0.8660)$]{
		\begin{tikzpicture}[scale=4]
			\draw[dashed] (0.5000,0.8660) -- (0,1) -- (-0.5000,0.8660);
			\draw[red,thick] (0,0) -- (0,1);
			\draw[blue,thick] (0,0) -- (0.5000,0.8660) -- (-0.5000,0.8660) -- cycle;
		\end{tikzpicture}
	}
	\subfloat[$(\geo{Q}_{8},3.1193,0.9730)$]{
		\begin{tikzpicture}[scale=4]
			\draw[dashed] (0,0) -- (0.3933,0.2424) -- (0.5000,0.6919) -- (0.3138,0.9495) -- (0,1) -- (-0.3138,0.9495) -- (-0.5000,0.6919) -- (-0.3933,0.2424) -- cycle;
			\draw[red,thick] (0,0) -- (0,1);
			\draw[blue,thick] (0,0) -- (0.3138,0.9495) -- (-0.3933,0.2424) -- (0.5000,0.6919) -- (-0.5000,0.6919) -- (0.3933,0.2424) -- (-0.3138,0.9495) -- cycle;
		\end{tikzpicture}
	}
	\subfloat[$(\geo{Q}_{16},3.1364,0.9942)$]{
		\begin{tikzpicture}[scale=4]
			\draw[dashed] (0,0) -- (0.2063,0.0604) -- (0.3738,0.1952) -- (0.4769,0.3838) -- (0.5000,0.5976) -- (0.4470,0.7665) -- (0.3333,0.9023) -- (0.1764,0.9843) -- (0,1) -- (-0.1764,0.9843) -- (-0.3333,0.9023) -- (-0.4470,0.7665) -- (-0.5000,0.5976) -- (-0.4769,0.3838) -- (-0.3738,0.1952) -- (-0.2063,0.0604) -- cycle;
			\draw[red,thick] (0,0) -- (0,1);
			\draw[blue,thick] (0,0) -- (0.1764,0.9843) -- (-0.2063,0.0604) -- (0.3333,0.9023) -- (-0.3738,0.1952) -- (0.4470,0.7665) -- (-0.4769,0.3838) -- (0.5000,0.5976) -- (-0.5000,0.5976) -- (0.4769,0.3838) -- (-0.4470,0.7665) -- (0.3738,0.1952) -- (-0.3333,0.9023) -- (0.2063,0.0604) -- (-0.1764,0.9843) -- cycle;
		\end{tikzpicture}
	}
	\subfloat[$(\geo{Q}_{32},3.1403,0.9987)$]{
		\begin{tikzpicture}[scale=4]
			\draw[dashed] (0,0) -- (0.1019,0.0149) -- (0.1989,0.0493) -- (0.2873,0.1020) -- (0.3637,0.1710) -- (0.4252,0.2535) -- (0.4695,0.3464) -- (0.4947,0.4462) -- (0.5000,0.5490) -- (0.4861,0.6413) -- (0.4544,0.7291) -- (0.4063,0.8091) -- (0.3434,0.8781) -- (0.2683,0.9335) -- (0.1838,0.9732) -- (0.0932,0.9956) -- (0,1) -- (-0.0932,0.9956) -- (-0.1838,0.9732) -- (-0.2683,0.9335) -- (-0.3434,0.8781) -- (-0.4063,0.8091) -- (-0.4544,0.7291) -- (-0.4861,0.6413) -- (-0.5000,0.5490) -- (-0.4947,0.4462) -- (-0.4695,0.3464) -- (-0.4252,0.2535) -- (-0.3637,0.1710) -- (-0.2873,0.1020) -- (-0.1989,0.0493) -- (-0.1019,0.0149) -- cycle;
			\draw[red,thick] (0,0) -- (0,1);
			\draw[blue,thick] (0,0) -- (0.0932,0.9956) -- (-0.1019,0.0149) -- (0.1838,0.9732) -- (-0.1989,0.0493) -- (0.2683,0.9335) -- (-0.2873,0.1020) -- (0.3434,0.8781) -- (-0.3637,0.1710) -- (0.4063,0.8091) -- (-0.4252,0.2535) -- (0.4544,0.7291) -- (-0.4695,0.3464) -- (0.4861,0.6413) -- (-0.4947,0.4462) -- (0.5000,0.5490) -- (-0.5000,0.5490) -- (0.4947,0.4462) -- (-0.4861,0.6413) -- (0.4695,0.3464) -- (-0.4544,0.7291) -- (0.4252,0.2535) -- (-0.4063,0.8091) -- (0.3637,0.1710) -- (-0.3434,0.8781) -- (0.2873,0.1020) -- (-0.2683,0.9335) -- (0.1989,0.0493) -- (-0.1838,0.9732) -- (0.1019,0.0149) -- (-0.0932,0.9956) -- cycle;
		\end{tikzpicture}
	}\\
	\subfloat[$(\geo{G}_8,3.1211,0.9776)$]{
		\begin{tikzpicture}[scale=4]
			\draw[dashed] (0,0) -- (0.2957,0.2043) -- (0.5000,0.5000) -- (0.4114,0.9114) -- (0,1) -- (-0.4114,0.9114) -- (-0.5000,0.5000) -- (-0.2957,0.2043) -- cycle;
			\draw[red,thick] (0,0) -- (0,1);
			\draw[blue,thick] (0,0) -- (0.4114,0.9114) -- (-0.5000,0.5000) -- (0.5000,0.5000) -- (-0.4114,0.9114) -- cycle;
			\draw (0.4114,0.9114) -- (-0.2957,0.2043);\draw (-0.4114,0.9114) -- (0.2957,0.2043);
		\end{tikzpicture}
	}
	\subfloat[$(\geo{G}_{16},3.1365,0.9950)$]{
		\begin{tikzpicture}[scale=4]
			\draw[dashed] (0,0) -- (0.1920,0.0578) -- (0.3473,0.1846) -- (0.4422,0.3613) -- (0.5000,0.5533) -- (0.4817,0.7439) -- (0.3598,0.8917) -- (0.1907,0.9817) -- (0,1) -- (-0.1907,0.9817) -- (-0.3598,0.8917) -- (-0.4817,0.7439) -- (-0.5000,0.5533) -- (-0.4422,0.3613) -- (-0.3473,0.1846) -- (-0.1920,0.0578) -- cycle;
			\draw[red,thick] (0,0) -- (0,1);
			\draw[blue,thick] (-0.5000,0.5533) -- (0.5000,0.5533);
			\draw[blue,thick] (0,0) -- (0.1907,0.9817) -- (-0.1920,0.0578) -- (0.3598,0.8917) -- (-0.3473,0.1846) -- (0.4817,0.7439) -- (-0.5000,0.5533);
			\draw[blue,thick] (0,0) -- (-0.1907,0.9817) -- (0.1920,0.0578) -- (-0.3598,0.8917) -- (0.3473,0.1846) -- (-0.4817,0.7439) -- (0.5000,0.5533);
			\draw (0.4817,0.7439) -- (-0.4422,0.3613);\draw (-0.4817,0.7439) -- (0.4422,0.3613);
		\end{tikzpicture}
	}
	\subfloat[$(\geo{G}_{32},3.1403,0.9988)$]{
		\begin{tikzpicture}[scale=4]
			\draw[dashed] (0,0) -- (0.0965,0.0143) -- (0.1884,0.0473) -- (0.2721,0.0975) -- (0.3444,0.1630) -- (0.4025,0.2414) -- (0.4527,0.3251) -- (0.4857,0.4170) -- (0.5000,0.5136) -- (0.4951,0.6121) -- (0.4711,0.7078) -- (0.4289,0.7970) -- (0.3627,0.8701) -- (0.2834,0.9289) -- (0.1943,0.9711) -- (0.0985,0.9951) -- (0,1) -- (-0.0985,0.9951) -- (-0.1943,0.9711) -- (-0.2834,0.9289) -- (-0.3627,0.8701) -- (-0.4289,0.7970) -- (-0.4711,0.7078) -- (-0.4951,0.6121) -- (-0.5000,0.5136) -- (-0.4857,0.4170) -- (-0.4527,0.3251) -- (-0.4025,0.2414) -- (-0.3444,0.1630) -- (-0.2721,0.0975) -- (-0.1884,0.0473) -- (-0.0965,0.0143) -- cycle;
			\draw[red,thick] (0,0) -- (0,1);
			\draw[blue,thick] (-0.5000,0.5136) -- (0.5000,0.5136);
			\draw[blue,thick] (0,0) -- (0.0985,0.9951) -- (-0.0965,0.0143) -- (0.1943,0.9711) -- (-0.1884,0.0473) -- (0.2834,0.9289) -- (-0.2721,0.0975) -- (0.3627,0.8701) -- (-0.3444,0.1630) -- (0.4289,0.7970) -- (-0.4527,0.3251) -- (0.4711,0.7078) -- (-0.4857,0.4170) -- (0.4951,0.6121) -- (-0.5000,0.5136);
			\draw[blue,thick] (0,0) -- (-0.0985,0.9951) -- (0.0965,0.0143) -- (-0.1943,0.9711) -- (0.1884,0.0473) -- (-0.2834,0.9289) -- (0.2721,0.0975) -- (-0.3627,0.8701) -- (0.3444,0.1630) -- (-0.4289,0.7970) -- (0.4527,0.3251) -- (-0.4711,0.7078) -- (0.4857,0.4170) -- (-0.4951,0.6121) -- (0.5000,0.5136);
			\draw (0.4289,0.7970) -- (-0.4025,0.2414);\draw (-0.4289,0.7970) -- (0.4025,0.2414);
		\end{tikzpicture}
	}
	\caption{Some families of convex small polygons $(\geo{P}_n,L(\geo{P}_n),W(\geo{P}_n))$}
	\label{figure:Pn}
\end{figure}

\section{Convex small polygons of large perimeter and large width}\label{sec:Sn}

We may ask if $\geo{J}_n$ is the best $n$-gon we can obtain with our approach presented in Section~\ref{sec:Cn} for any $n = 2^s$ with integer $s\ge 2$. This section answers this question.

Let $n=2^s \ge 4$. For an odd integer $m=3,5,\ldots, n-1$ and for a composition $(c_0, c_1,\ldots, c_{\frac{m-1}{2}})$ of $n/2$ of size $\lceil m/2 \rceil$, consider the function
\begin{equation}
	\label{eq:violation}
	d(c_0, c_1,\ldots, c_{\frac{m-1}{2}}) := \left| \sum_{j=1}^{\frac{m-1}{2}} (-1)^{j-1} \sin \left(\sum_{i=0}^{j-1} c_i \frac{\pi}{n}\right) + \frac{(-1)^{{\frac{m-1}{2}}}}{2} \right|
\end{equation}
that measures the violation of the cycle constraint~\eqref{eq:x4} by the solution
\[
(\alpha_0, \alpha_1,\ldots, \alpha_{\frac{m-1}{2}}) = (\pi/n, \pi/n,\ldots,\pi/n)
\]
or~\eqref{eq:x4d} by $\delta_0 = 0$. Let $(c_0^*, c_1^*,\ldots, c_{\frac{m-1}{2}}^*)$ denote the composition that minimizes the constraint violation $d$ in~\eqref{eq:violation} among all $\binom{n/2-1}{\lfloor m/2 \rfloor}$ compositions $(c_0, c_1,\ldots, c_{\frac{m-1}{2}})$ of $n/2$ of size $\lceil m/2 \rceil$ and let $\geo{S}_{m,n}$ denote the corresponding $n$-gon obtained by solving~\eqref{eq:x4d}. We observed that the optimal composition $(c_0^*, c_1^*,\ldots, c_{\frac{m-1}{2}}^*)$ simultaneously maximizes the perimeter $L(\geo{P}_n)$ and the width $W(\geo{P}_n)$ in~\eqref{eq:LWd}. If $m=n-1$, we have $\geo{S}_{n-1,n} \equiv \geo{Q}_n$ since there exists only one composition of $n/2$ of size $n/2$. If $m=3$, we can show that $\geo{S}_{3,n} \equiv \geo{T}_n$. If $m=n-3$, we remarked that $\geo{S}_{n-3,n} \equiv \geo{G}_n$.

Now, let $\geo{S}_n$ denote the best $n$-gon among all $n$-gons $\geo{S}_{m,n}$. The composition of~$\geo{S}_n$ minimizes the constraint violation $d$~\eqref{eq:violation} among all $2^{n/2-1}-1$ compositions of $n/2$. We computed $\geo{S}_n$, illustrated in Figure~\ref{figure:Sn}, for $n=4,8,16,32,64$ and its characteristics are summarized in Table~\ref{table:Sn}. We note that the composition $(c_0^*,c_1^*,\ldots,c_{\frac{m-1}{2}}^*)$ of $\geo{S}_n$ verifies
\[
\sum_{k=0}^{\frac{m-1}{2}} (-1)^k c_k^* = 0.
\]
We also remark that $\geo{S}_n \equiv \geo{J}_n$ for $n =4,8,16$, i.e., $\geo{J}_n$ is the best $n$-gon we can obtain with our approach when $n\le 16$. That is not the case when $n\ge 32$.

\begin{table}[h]
	\footnotesize
	\centering
	\caption{Characteristics of $\geo{S}_n$}
	\label{table:Sn}
		\begin{tabular}{@{}rrll@{}}
			\toprule
			$n$ & $m$ & $(c_0^*,c_1^*,\ldots,c_{\frac{m-1}{2}}^*)$ & $\delta_0$ \\
			\midrule
			4	&	3	&	$(1,1)$	&	$-0.2617993877991494$	\\
			8	&	5	&	$(1,2,1)$	&	$\hphantom{-}0.3133195779201635\cdot 10^{-1}$	\\
			16	&	11	&	$(1,2,1,1,2,1)$	&	$\hphantom{-}0.1583040448658572\cdot 10^{-2}$	\\
			32	&	21	&	$(1,1,2,2,1,1,1,1,1,3,2)$	&	$-0.6544777934135723\cdot 10^{-6}$	\\
			64	&	27	&	$(1,6,5,1,2,1,1,3,2,1,3,1,2,3)$	&	$\hphantom{-}0.9512026924844603\cdot 10^{-11}$	\\
			\bottomrule
		\end{tabular}
\end{table}

The perimeters $L(\geo{S}_n)$ of $\geo{S}_n$ are displayed in Table~\ref{table:L(Sn)} and the widths $W(\geo{S}_n)$ in Table~\ref{table:W(Sn)}. We see that, when $n\ge 32$, the $n$-gon $\geo{S}_n$ outperforms $\geo{J}_n$. For $n=64$, we note
\[
\begin{aligned}
	L_{64}^* - L(\geo{S}_{64}) &< \ub{L}_{64} - L(\geo{S}_{64}) &&= 3.56\ldots \times 10^{-23}\\
	&< \ub{L}_{64} - L(\geo{J}_{64}) &&= 9.05\ldots \times 10^{-14},\\
	W_{64}^* - W(\geo{S}_{64}) &< \ub{W}_{64} - W(\geo{S}_{64}) &&= 1.17\ldots \times 10^{-13}\\
	&< \ub{W}_{64} - W(\geo{J}_{64}) &&= 5.89\ldots \times 10^{-9}.
\end{aligned}
\]

\begin{table}[h]
	\footnotesize
	\centering
	\caption{Perimeters of $\geo{S}_n$}
	\label{table:L(Sn)}
		\begin{tabular}{@{}rllll@{}}
			\toprule
			$n$ & $L(\geo{R}_n)$ & $L(\geo{J}_n)$ & $L(\geo{S}_n)$ & $\ub{L}_n$ \\
			\midrule
			4	&	2.83	&	3.0353	&	3.0353	&	3.0615	\\
			8	&	3.06	&	3.121062	&	3.121062	&	3.121445	\\
			16	&	3.1214	&	3.13654751	&	3.13654751	&	3.13654849	\\
			32	&	3.1365	&	3.140331156355	&	3.14033115695458	&	3.14033115695475	\\
			64	&	3.1403	&	3.14127725093268	&	3.141277250932772868061984	&	3.141277250932772868062019	\\
			\bottomrule
		\end{tabular}
\end{table}

\begin{table}[h]
	\footnotesize
	\centering
	\caption{Widths of $\geo{S}_n$}
	\label{table:W(Sn)}
		\begin{tabular}{@{}rllll@{}}
			\toprule
			$n$ & $W(\geo{R}_n)$ & $W(\geo{J}_n)$ & $W(\geo{S}_n)$ & $\ub{W}_n$ \\
			\midrule
			4	&	0.71	&	0.87	&	0.87	&	0.92	\\
			8	&	0.9239	&	0.9776	&	0.9776	&	0.9808	\\
			16	&	0.9808	&	0.995107	&	0.995107	&	0.995185	\\
			32	&	0.9952	&	0.99879450	&	0.9987954401	&	0.9987954562	\\
			64	&	0.9988	&	0.9996988128	&	0.99969881869609	&	0.99969881869620	\\
			\bottomrule
		\end{tabular}
\end{table}

\begin{figure}[H]
	\centering
	\subfloat[$(\geo{S}_4,3.035276,0.866025)$]{
		\begin{tikzpicture}[scale=5]
			\draw[dashed] (0.5000,0.8660) -- (0,1) -- (-0.5000,0.8660);
			\draw[blue,thick] (0,0) -- (0.5000,0.8660) -- (-0.5000,0.8660) -- cycle;
			\draw[red,thick] (0,0) -- (0,1);
		\end{tikzpicture}
	}
	\subfloat[$(\geo{S}_8,3.121062,0.977609)$]{
		\begin{tikzpicture}[scale=5]
			\draw[dashed] (0,0) -- (0.2957,0.2043) -- (0.5000,0.5000) -- (0.4114,0.9114) -- (0,1) -- (-0.4114,0.9114) -- (-0.5000,0.5000) -- (-0.2957,0.2043) -- cycle;
			\draw[blue,thick] (0,0) -- (0.4114,0.9114) -- (-0.5000,0.5000) -- (0.5000,0.5000) -- (-0.4114,0.9114) -- cycle;
			\draw[red,thick] (0,0) -- (0,1);
			\draw (0.4114,0.9114) -- (-0.2957,0.2043);\draw (-0.4114,0.9114) -- (0.2957,0.2043);
		\end{tikzpicture}
	}
	\subfloat[$(\geo{S}_{16},3.136548,0.995107)$]{
		\begin{tikzpicture}[scale=5]
			\draw[dashed] (0,0) -- (0.1860,0.0566) -- (0.3576,0.1481) -- (0.4811,0.2983) -- (0.5000,0.4919) -- (0.4428,0.6810) -- (0.3495,0.8552) -- (0.1966,0.9805) -- (0,1) -- (-0.1966,0.9805) -- (-0.3495,0.8552) -- (-0.4428,0.6810) -- (-0.5000,0.4919) -- (-0.4811,0.2983) -- (-0.3576,0.1481) -- (-0.1860,0.0566) -- cycle;
			\draw[red,thick] (0,0)--(0,1);
			\draw[blue,thick] (0.5000,0.4919)--(-0.5000,0.4919);
			\draw[blue,thick] (0,0) -- (0.1966,0.9805) -- (-0.3576,0.1481) -- (0.3495,0.8552) -- (-0.4847,0.3006) -- (0.5000,0.4919);
			\draw[blue,thick] (0,0) -- (-0.1966,0.9805) -- (0.3576,0.1481) -- (-0.3495,0.8552) -- (0.4847,0.3006) -- (-0.5000,0.4919);
			\draw (0.1966,0.9805) -- (-0.1860,0.0566);\draw (-0.1966,0.9805) -- (0.1860,0.0566);
			\draw (-0.4847,0.3006) -- (0.4428,0.6810);\draw (0.4847,0.3006) -- (-0.4428,0.6810);
		\end{tikzpicture}
	}\\
	\subfloat[$(\geo{S}_{32},3.140331,0.998795)$]{
		\begin{tikzpicture}[scale=7]
			\draw[dashed] (0,0) -- (0.0971,0.0144) -- (0.1858,0.0564) -- (0.2700,0.1068) -- (0.3427,0.1727) -- (0.4011,0.2515) -- (0.4431,0.3403) -- (0.4762,0.4327) -- (0.5000,0.5278) -- (0.4952,0.6259) -- (0.4808,0.7229) -- (0.4303,0.8071) -- (0.3644,0.8798) -- (0.2856,0.9383) -- (0.1932,0.9713) -- (0.0980,0.9952) -- (0,1) -- (-0.0980,0.9952) -- (-0.1932,0.9713) -- (-0.2856,0.9383) -- (-0.3644,0.8798) -- (-0.4303,0.8071) -- (-0.4808,0.7229) -- (-0.4952,0.6259) -- (-0.5000,0.5278) -- (-0.4762,0.4327) -- (-0.4431,0.3403) -- (-0.4011,0.2515) -- (-0.3427,0.1727) -- (-0.2700,0.1068) -- (-0.1858,0.0564) -- (-0.0971,0.0144) -- cycle;
			\draw[red,thick] (0,0)--(0,1);
			\draw[blue,thick] (-0.5000,0.5278)--(0.5000,0.5278);
			\draw[blue,thick] (0,0) -- (0.0980,0.9952) -- (-0.0971,0.0144) -- (0.2856,0.9383) -- (-0.2700,0.1068) -- (0.3644,0.8798) -- (-0.3427,0.1727) -- (0.4303,0.8071) -- (-0.4011,0.2515) -- (0.4808,0.7229) -- (-0.5000,0.5278);
			\draw[blue,thick] (0,0) -- (-0.0980,0.9952) -- (0.0971,0.0144) -- (-0.2856,0.9383) -- (0.2700,0.1068) -- (-0.3644,0.8798) -- (0.3427,0.1727) -- (-0.4303,0.8071) -- (0.4011,0.2515) -- (-0.4808,0.7229) -- (0.5000,0.5278);
			\draw (0.0971,0.0144) -- (-0.1932,0.9713);\draw (-0.0971,0.0144) -- (0.1932,0.9713);
			\draw (-0.2856,0.9383) -- (0.1858,0.0564);\draw (0.2856,0.9383) -- (-0.1858,0.0564);
			\draw (0.4808,0.7229) -- (-0.4431,0.3403);\draw (-0.4808,0.7229) -- (0.4431,0.3403);
			\draw (0.4808,0.7229) -- (-0.4762,0.4327);\draw (-0.4808,0.7229) -- (0.4762,0.4327);
			\draw (-0.5000,0.5278) -- (0.4952,0.6259);\draw (0.5000,0.5278) -- (-0.4952,0.6259);
		\end{tikzpicture}
	}
	\subfloat[$(\geo{S}_{64},3.141277,0.999699)$]{
		\begin{tikzpicture}[scale=7]
			\draw[dashed] (0,0) -- (0.0489,0.0036) -- (0.0977,0.0096) -- (0.1460,0.0180) -- (0.1939,0.0288) -- (0.2412,0.0419) -- (0.2878,0.0573) -- (0.3280,0.0855) -- (0.3635,0.1194) -- (0.3956,0.1565) -- (0.4258,0.1952) -- (0.4540,0.2353) -- (0.4761,0.2792) -- (0.4892,0.3265) -- (0.4952,0.3752) -- (0.4988,0.4242) -- (0.5000,0.4732) -- (0.4916,0.5216) -- (0.4809,0.5695) -- (0.4655,0.6161) -- (0.4478,0.6619) -- (0.4279,0.7067) -- (0.4037,0.7494) -- (0.3775,0.7909) -- (0.3436,0.8265) -- (0.3064,0.8585) -- (0.2677,0.8887) -- (0.2263,0.9150) -- (0.1836,0.9392) -- (0.1397,0.9612) -- (0.0949,0.9811) -- (0.0491,0.9988) -- (0,1) -- (-0.0491,0.9988) -- (-0.0949,0.9811) -- (-0.1397,0.9612) -- (-0.1836,0.9392) -- (-0.2263,0.9150) -- (-0.2677,0.8887) -- (-0.3064,0.8585) -- (-0.3436,0.8265) -- (-0.3775,0.7909) -- (-0.4037,0.7494) -- (-0.4279,0.7067) -- (-0.4478,0.6619) -- (-0.4655,0.6161) -- (-0.4809,0.5695) -- (-0.4916,0.5216) -- (-0.5000,0.4732) -- (-0.4988,0.4242) -- (-0.4952,0.3752) -- (-0.4892,0.3265) -- (-0.4761,0.2792) -- (-0.4540,0.2353) -- (-0.4258,0.1952) -- (-0.3956,0.1565) -- (-0.3635,0.1194) -- (-0.3280,0.0855) -- (-0.2878,0.0573) -- (-0.2412,0.0419) -- (-0.1939,0.0288) -- (-0.1460,0.0180) -- (-0.0977,0.0096) -- (-0.0489,0.0036) -- cycle;
			\draw[red,thick] (0,0)--(0,1);
			\draw[blue,thick] (0.5000,0.4732)--(-0.5000,0.4732);
			\draw[blue,thick] (0,0) -- (0.0491,0.9988) -- (-0.2878,0.0573) -- (0.2677,0.8887) -- (-0.3280,0.0855) -- (0.3436,0.8265) -- (-0.3635,0.1194) -- (0.3775,0.7909) -- (-0.4540,0.2353) -- (0.4279,0.7067) -- (-0.4761,0.2792) -- (0.4809,0.5695) -- (-0.4892,0.3265) -- (0.5000,0.4732);
			\draw[blue,thick] (0,0) -- (-0.0491,0.9988) -- (0.2878,0.0573) -- (-0.2677,0.8887) -- (0.3280,0.0855) -- (-0.3436,0.8265) -- (0.3635,0.1194) -- (-0.3775,0.7909) -- (0.4540,0.2353) -- (-0.4279,0.7067) -- (0.4761,0.2792) -- (-0.4809,0.5695) -- (0.4892,0.3265) -- (-0.5000,0.4732);
			\draw (0.0491,0.9988) -- (-0.0489,0.0036);\draw (-0.0491,0.9988) -- (0.0489,0.0036);
			\draw (0.0491,0.9988) -- (-0.0977,0.0096);\draw (-0.0491,0.9988) -- (0.0977,0.0096);
			\draw (0.0491,0.9988) -- (-0.1460,0.0180);\draw (-0.0491,0.9988) -- (0.1460,0.0180);
			\draw (0.0491,0.9988) -- (-0.1939,0.0288);\draw (-0.0491,0.9988) -- (0.1939,0.0288);
			\draw (0.0491,0.9988) -- (-0.2412,0.0419);\draw (-0.0491,0.9988) -- (0.2412,0.0419);
			\draw (-0.2878,0.0573) -- (0.0949,0.9811);\draw (0.2878,0.0573) -- (-0.0949,0.9811);
			\draw (-0.2878,0.0573) -- (0.1397,0.9612);\draw (0.2878,0.0573) -- (-0.1397,0.9612);
			\draw (-0.2878,0.0573) -- (0.1836,0.9392);\draw (0.2878,0.05730) -- (-0.1836,0.9392);
			\draw (-0.2878,0.0573) -- (0.2263,0.9150);\draw (0.2878,0.0573) -- (-0.2263,0.9150);
			\draw (-0.3280,0.0855) -- (0.3064,0.8585);\draw (0.3280,0.0855) -- (-0.3064,0.8585);
			\draw (0.3775,0.7909) -- (-0.3956,0.1565);\draw (-0.3775,0.7909) -- (0.3956,0.1565);
			\draw (0.3775,0.7909) -- (-0.4258,0.1952);\draw (-0.3775,0.7909) -- (0.4258,0.1952);
			\draw (-0.4540,0.2353) -- (0.4037,0.7494);\draw (0.4540,0.2353) -- (-0.4037,0.7494);
			\draw (-0.4761,0.2792) -- (0.4478,0.6619);\draw (0.4761,0.2792) -- (-0.4478,0.6619);
			\draw (-0.4761,0.2792) -- (0.4655,0.6161);\draw (0.4761,0.2792) -- (-0.4655,0.6161);
			\draw (-0.4892,0.3265) -- (0.4916,0.5216);\draw (0.4892,0.3265) -- (-0.4916,0.5216);
			\draw (0.5000,0.4732) -- (-0.4952,0.3752);\draw (-0.5000,0.4732) -- (0.4952,0.3752);
			\draw (0.5000,0.4732) -- (-0.4988,0.4242);\draw (-0.5000,0.4732) -- (0.4988,0.4242);
		\end{tikzpicture}
	}
	\caption{Polygons $(\geo{S}_n,L(\geo{S}_n),W(\geo{S}_n))$: (a) Tetragon $\geo{S}_4$; (b) Octagon $\geo{S}_8$; (c) Hexadecagon $\geo{S}_{16}$; (d) Triacontadigon $\geo{S}_{32}$; (e) Hexacontatetragon $\geo{S}_{64}$}
	\label{figure:Sn}
\end{figure}

\section{Maximizing the perimeter}\label{sec:maxperi}
For $n=2^s$ with $s\ge 3$, we can improve the perimeter of $\geo{Q}_n$, $\geo{B}_n$, $\geo{G}_n$, $\geo{J}_n$ or $\geo{S}_n$, by adjusting the angles~$\alpha_k$ from our parametrization of Section~\ref{sec:Cn} to maximize the perimeter $L(\geo{P}_n)$ in~\eqref{eq:LW:L}, creating a polygon $\geo{P}_n^*$ with larger perimeter. Thus, $L(\geo{P}_n^*)$ is the optimal value of the following optimization problem:
\begin{subequations}\label{eq:ngon:Sn:L}
	\begin{align}
		L(\geo{P}_n^*) = \max_{\rv{\alpha}} \quad & \sum_{k=0}^{\frac{m-1}{2}} 4c_k\sin \frac{\alpha_k}{2}\\
		\subj \quad & \sum_{k=0}^{\frac{m-1}{2}} c_k \alpha_k = \pi/2,\\
		& \sum_{j=1}^{\frac{m-1}{2}} (-1)^{j-1} \sin \left(\sum_{i=0}^{j-1} c_i \alpha_i\right) = (-1)^{{\frac{m+1}{2}}}/2,\\
		& 0 \le c_0 \alpha_0 \le \pi/6,\\
		& 0 \le c_k \alpha_k \le \pi/3 \quad \forall k = 1,2, \ldots, {\frac{m-1}{2}},
	\end{align}
\end{subequations}
where $(c_0,c_1,\ldots,c_{\frac{m-1}{2}})$ is the composition of $n/2$ corresponding to $\geo{P}_n$.

Note that $L(\geo{B}_8^*) = L_8^*$. Then Mossinghoff asked if $L(\geo{B}_{16}^*) = L_{16}^*$ and if the maximal perimeter when $n = 2^s$ is always achieved by a polygon with the same diameter graph as $\geo{B}_n$. Numerical results in Table~\ref{table:L(Pn*)} show that both conjectures might be false. Indeed, for each $n$, we have $L(\geo{B}_n^*) < L(\geo{C}_n) < L(\geo{C}_n^*)$.

Problem~\eqref{eq:ngon:Sn:L} was solved in Julia using JuMP~\cite{dunning2017} with Couenne~0.5.8. Julia codes have been made available in OPTIGON. The solver Couenne~\cite{belotti2009} is a branch-and-bound algorithm that aims at finding global optima of nonconvex mixed-integer nonlinear optimization problems. This means that, under the composition of $n/2$ describing $\geo{P}_n$, $L(\geo{P}_n^*)$ is the longest perimeter.

Table~\ref{table:L(Pn*)} shows the optimal values $L(\geo{Q}_n^*)$, $L(\geo{B}_n^*)$, $L(\geo{G}_n^*)$, $L(\geo{C}_n^*)$, $L(\geo{J}_n^*)$ along with the upper bounds~$\ub{L}_n$ for $n=8,16,32,64$. The optimal values $L(\geo{S}_n^*)$ are given in Table~\ref{table:L(Sn*)} with $L(\geo{S}_n)$ and~$\ub{L}_n$ for $n=8,16,32$. We could not compute $L(\geo{S}_{64}^*)$ due to the limit of numerical computations. We report in Table~\ref{table:g(Pn*)} the fractions $\gamma (\geo{P}_n^*) := \frac{L(\geo{P}_n^*) - L(\geo{P}_n)}{\ub{L}_n-L(\geo{P}_n)}$ of the length of the interval $[L(\geo{P}_n), \ub{L}_n]$ where $L(\geo{P}_n^*)$ lies. The results support the following key points:
\begin{enumerate}
	\item For $n=8$, $L(\geo{B}_8^*) = L(\geo{G}_8^*) = L(\geo{J}_8^*) = L(\geo{S}_8^*)  = L_8^*$ according to Theorem~\ref{thm:perimeter}~\cite{griffiths1975,audet2007a}.
	\item For each $n \ge 16$, $L(\geo{B}_n^*) < L(\geo{C}_n) < L(\geo{C}_n^*)$, i.e., $\geo{B}_n^*$ is a suboptimal solution.
	\item For each $n \ge 32$, $L(\geo{J}_n^*) < L(\geo{S}_n)$, i.e., $\geo{J}_n^*$ is a suboptimal solution.
	\item For each $\geo{P}_n^*$, the fraction $\gamma (\geo{P}_n^*)$ in Table~\ref{table:g(Pn*)} appears to approach to the same scalar $\gamma^* \in (0,1)$ as $n$ increases, i.e., $\ub{L}_n - L(\geo{P}_n^*) = O(\ub{L}_n - L(\geo{P}_n))$. We suspect $\gamma^* = 1-\frac{8}{\pi^2} =0.1894\ldots$ since for large $n$,
	\[
	\ub{L}_n - L(\geo{Q}_n^*) \approx \ub{L}_n - L(\geo{T}_n) \sim \frac{\pi^3}{4n^4} \sim \frac{8}{\pi^2} (\ub{L}_n - L(\geo{Q}_n)).
	\]
\end{enumerate}
Assuming the existence of an axis of symmetry, it follows from Proposition~\ref{thm:SnL}, $\geo{S}_n^*$ has the largest perimeter among all convex small $n$-gons. Recently, Mulansky and Potschka~\cite{mulansky2022} claimed that $L(\geo{S}_{16}^*) = L_{16}^*$ and $L(\geo{S}_{32}^*) = L_{32}^*$. We then asked if $L(\geo{S}_n^*) = L_n^*$ when $n$ is a power of $2$.

\begin{proposition}
	\label{thm:SnL}
	Let $n=8,16,32$. Under the assumption of symmetry, $\geo{S}_n^*$ has the largest perimeter among all convex small $n$-gons.
\end{proposition}
\begin{proof}
	For $n=8,16,32$, solving Problem~\eqref{eq:ngon:Sn:L} with the solver Couenne for each composition of $n/2$, we obtained that $\geo{S}_n^*$ has the largest perimeter within the limit of the numerical computations.
\end{proof}

\begin{table}[H]
	\footnotesize
	\centering
	\caption{Perimeters of $\geo{Q}_n^*$, $\geo{B}_n^*$, $\geo{G}_n^*$, $\geo{C}_n^*$, and $\geo{J}_n^*$}
	\label{table:L(Pn*)}
	\resizebox{\linewidth}{!}{
		\begin{tabular}{@{}rllllll@{}}
			\toprule
			$n$ & $L(\geo{Q}_n^*)$ & $L(\geo{B}_n^*)$ & $L(\geo{G}_n^*)$ & $L(\geo{C}_n^*)$ & $L(\geo{J}_n^*)$ & $\ub{L}_n$ \\
			\midrule
			8	&	3.1196	&	3.121147	&	3.121147	&	--	&	3.121147	&	3.121445	\\
			16	&	3.136431	&	3.13654396	&	3.13654225	&	3.13654772	&	3.13654772	&	3.13654849	\\
			32	&	3.140324	&	3.14033109	&	3.14033107	&	3.1403311541	&	3.140331156472	&	3.140331156955	\\
			64	&	3.14127679	&	3.1412772498	&	3.1412772494	&	3.141277250922	&	3.1412772509326996	&	3.1412772509327729	\\
			\bottomrule
		\end{tabular}
	}
\end{table}

\begin{table}[h]
	\footnotesize
	\centering
	\caption{Perimeters of $\geo{S}_n^*$}
	\label{table:L(Sn*)}
		\begin{tabular}{@{}rlll@{}}
			\toprule
			$n$ & $L(\geo{S}_n)$ & $L(\geo{S}_n^*)$ & $\ub{L}_n$	\\
			\midrule
			8	&	3.121062	&	3.121147	&	3.121445	\\
			16	&	3.13654751	&	3.13654772	&	3.13654849	\\
			32	&	3.14033115695458	&	3.14033115695462	&	3.14033115695475	\\
			\bottomrule
		\end{tabular}
\end{table}

\begin{table}[H]
	\footnotesize
	\centering
	\caption{Fractions $\gamma (\geo{P}_n^*) := \frac{L(\geo{P}_n^*) - L(\geo{P}_n)}{\ub{L}_n-L(\geo{P}_n)}$}
	\label{table:g(Pn*)}
		\begin{tabular}{@{}rllllll@{}}
			\toprule
			$n$ & $\gamma(\geo{Q}_n^*)$ & $\gamma(\geo{B}_n^*)$ & $\gamma(\geo{G}_n^*)$ & $\gamma(\geo{C}_n^*)$ & $\gamma(\geo{J}_n^*)$	&	$\gamma(\geo{S}_n^*)$ \\
			\midrule
			8	&	0.1231	&	0.2219	&	0.2219	&	--	&	0.2219	&	0.2219	\\
			16	&	0.1724	&	0.2077	&	0.2157	&	0.2122	&	0.2122	&	0.2122	\\
			32	&	0.1851	&	0.1945	&	0.1948	&	0.1947	&	0.1947	&	0.2026	\\
			64	&	0.1884	&	0.1907	&	0.1911	&	0.1908	&	0.1909	&	--	\\
			\bottomrule
		\end{tabular}
\end{table}

\begin{table}[H]
	\footnotesize
	\centering
	\caption{Angles $\alpha_0^*, \alpha_1^*, \ldots, \alpha_{\frac{m-1}{2}}^*$ of $\geo{S}_n^*$}
	\label{table:Sn:angles}
		\begin{tabular}{@{}ll|lllllll@{}}
			\toprule
			$n$ & $\pi/n$ & $i$ & $\alpha_{6i}^*$ & $\alpha_{6i+1}^*$ & $\alpha_{6i+2}^*$ & $\alpha_{6i+3}^*$ & $\alpha_{6i+4}^*$ & $\alpha_{6i+5}^*$ \\
			\midrule
			8	&	0.392699	&	0	&	0.435281	&	0.368535	&	0.398447	&&&\\
			16	&	0.196350	&	0	&	0.198316	&	0.194503	&	0.197746	&	0.194994	&	0.197164	&	0.196406	\\
			32	&	0.0981748	&	0	&	0.0981739	&	0.0981755	&	0.0981739	&	0.0981754	&	0.0981741	&	0.0981754	\\
			&& 1	&	0.0981742 & 0.0981752 & 0.0981743 & 0.0981751 & 0.0981748 &\\
			\bottomrule
		\end{tabular}
\end{table}

The optimal angles $\alpha_k^*$ that produce $\geo{S}_n^*$ appear in Table~\ref{table:Sn:angles}. They exhibit a pattern of damped oscillation, converging in an alterning manner to a mean value around $\pi/n$. We remark that
\[
W(\geo{S}_n^*) = \min_k \cos (\alpha_k ^*/2) < W(\geo{S}_n)
\]
for each $n$. From Theorem~\ref{thm:perimeter}, $W(\geo{S}_n) = W_n^*$ when $n=4$ and $n=8$~\cite{bezdek2000,audet2013}. However, from Proposition~\ref{thm:S16W}, $W(\geo{S}_{16}) < W_{16}^*$ as we construct a convex small $16$-gon $\geo{S}_{16}^{\star}$ whose width exceeds that of~$\geo{S}_{16}$. We ask if $W(\geo{S}_{16}^{\star}) = W_{16}^*$.

\begin{proposition}
	\label{thm:S16W}
	There exists a convex small $16$-gon $\geo{S}_{16}^{\star}$ such that
	\[
	W(\geo{S}_{16}^{\star}) = 0.995132\ldots > W(\geo{S}_{16}).
	\]
\end{proposition}
\begin{proof}
	Consider a convex small $16$-gon $\geo{P}_{16}$ having the same diameter graph as $\geo{S}_{16}$. Thus, the $16$-gon~$\geo{P}_{16}$ is described by the composition $(c_0,c_1,c_2,c_4,c_5)=(1,2,1,1,2,1)$ of $8$. Suppose
	\[
	\alpha_k =
	\begin{cases}
		\beta & \text{if $k=1$,}\\
		\alpha & \text{otherwise.}
	\end{cases}
	\]
	Then, \eqref{eq:condition} and \eqref{eq:x4} become
	\[
	\begin{cases}
		6\alpha + 2\beta = \pi/2,\\
		\sin (\alpha) - \sin (\alpha+2\beta) + \sin (2\alpha+2\beta)-\sin (3\alpha+2\beta) + \sin (5\alpha+2\beta) = 1/2,
	\end{cases}
	\]
	respectively. This system of equations has a solution
	\[
	(\hat{\alpha},\hat{\beta}) = (0.197421\ldots,0.193134\ldots)
	\]
	Let $\geo{S}_{16}^{\star}$ denote the $16$-gon obtained by setting $(\alpha,\beta)=(\hat{\alpha},\hat{\beta})$. The width of $\geo{S}_{16}^{\star}$ is
	\[
	W(\geo{S}_{16}^{\star}) = \cos(\hat{\alpha}/2) = 0.995132\ldots > W(\geo{S}_{16}) = 0.995107\ldots
	\]
	and its perimeter is
	\[
	L(\geo{S}_{16}^{\star}) = 24\sin(\hat{\alpha}/2) + 8\sin(\hat{\beta}/2) = 3.136542\ldots < L(\geo{S}_{16}).
	\]
\end{proof}

\section{Conclusion}\label{sec:conclusion}
We provided tigther lower bounds on the maximal perimeter and the maximal width of convex small $n$-gons when $n$ is a power of $2$. Under the assumption of symmetry, we described a convex small $n$-gon as a composition of $n/2$ and we expressed both its perimeter and its width as functions of a single variable.
The convex small $n$-gon of large perimeter and large width is then obtained by selecting the composition of $n/2$ that minimizes the violation of a cycle constraint by a particular solution. The family of the $n$-gons constructed outperforms any other family found in the literature. From our results, it appears that Mossinghoff conjecture on the diameter graph of a convex small $n$-gon of maximal perimeter is not true when $n\ge 16$.

\section*{Acknowledgements}
The author thanks Charles Audet, Professor at Polytechnique Montr\'{e}al, and Kevin Hare, Professor at University of Waterloo, for helpful discussions on extremal small polygons.

\bibliographystyle{ieeetr}
\bibliography{../../research}
\end{document}